%% file: GMP.tex
\documentclass[12pt]{amsart}
\usepackage{graphicx}

\usepackage{amscd,amsmath,amssymb}
\newtheorem{theorem}{Theorem}[section]
\newtheorem{proposition}[theorem]{Proposition}
\newtheorem{lemma}[theorem]{Lemma}

\newtheorem{corollary}[theorem]{Corollary}
\theoremstyle{definition}
\newtheorem{definition}[theorem]{Definition}
\newtheorem{example}[theorem]{Example}

\newcommand{\aaa}{\mbox{$\alpha$}}

\newcommand{\bbb}{\mbox{$\beta$}}

\newcommand{\bdd}{\mbox{$\partial$}}

\theoremstyle{remark}
\newtheorem{remark}[theorem]{Remark}

\numberwithin{equation}{section}

\input auxxx.tex

\begin{document}

\title[Heegaard splittings and closed orbits of gradient flows]
{Morse-Novikov theory, Heegaard splittings and 
 closed orbits of gradient flows} 

\date{\today}

\begin{abstract}
The works of Donaldson \cite{donaldson} and Mark \cite{mark} 
make the structure of the Seiberg-Witten invariant of 3-manifolds clear. 
It corresponds to certain torsion type invariants counting 
flow lines and closed orbits of a gradient flow 
of a circle-valued Morse map on a 3-manifold.  
We study these invariants using the Morse-Novikov theory and 
Heegaard splitting for sutured manifolds, and 
make detailed computations for knot complements. 
\end{abstract}

\author{Hiroshi Goda}
\address{Department of Mathematics,
Tokyo University of Agriculture and Technology,
2-24-16 Naka-cho, Koganei,
Tokyo 184-8588, Japan}
\email{goda@cc.tuat.ac.jp}
\author{Hiroshi Matsuda}
\address{Department of Mathematics, 
Graduate school of Science,
Hiroshima University,
Hiroshima 739-8526, Japan}
\email{matsuda@math.sci.hiroshima-u.ac.jp}
\author{Andrei Pajitnov}
\address{Laboratoire de Math\'ematiques Jean-Leray UMR 6629,
Universit\'e de Nantes, Facult\'e des Sciences,   
2, rue de la Houssini\`ere,  
44072, Nantes, Cedex,  
France}
\email{pajitnov@math.univ-nantes.fr}

\thanks{
The first and second authors are partially supported
by Grant-in-Aid for Scientific Research,
(No.21540071 and No.20740041),
Ministry of Education, Science, 
Sports and Technology, Japan.
}

\maketitle

%--------------Section1.tex-----'07/09/18-----

\section{Introduction}
\label{s:intro}

Let $K\subset S^3$ be an oriented knot, 
put $C_{K}=S^3-K$. 
The canonical cohomology class $\xi\in H^{1}(C_{K})=[C_{K},S^{1}]$ 
can be represented by a Morse map $f: C_{K} \to S^{1}.$ 
In this paper we study the dynamics of the gradient flow of $f$. 

Milnor pointed out in \cite{milnor2} 
a relationship between the Reidemeister torsion and 
dynamical  zeta functions. 
His theorem applies to fibred knots, that is to the case 
when we can choose the map $f$ without critical points. 
The theorem implies in particular that 
the Alexander polynomial of any fibred knot in $S^{3}$ 
is essentially the same as
the Lefschetz zeta function of the monodromy 
map of the fibration $f$. 
The periodic points of the monodromy map correspond 
to the closed orbits of the gradient flow of the 
fibration $C_{K}\to S^{1}$; 
thus Milnor's theorem establishes a relation 
between the dynamics of this gradient flow and 
and the Alexander polynomial of the knot. 

When the knot $K$ is not fibred, the Morse map $f$ 
necessarily has critical points. 
The Milnor's formula is no more valid, however it can 
be generalized to this case at the cost of adding a correction term, 
as it was discovered by Hutchings and Lee 
(\cite{hulee1}, \cite{hulee2}). 
This correction term is essentially the torsion 
of the Novikov complex associated with 
the circle-valued Morse map $f$ 
(see \cite{novikov}, \cite{pajitnov06}). 
This complex is an analog of the Morse complex 
for the circle-valued case, 
and is obtained through counting the flow lines 
of the gradient joining the critical points 
of the map. 

The torsion of the Novikov complex and the Lefschetz zeta 
function are in general very difficult to compute 
due to the complexity of the transversal gradient flows 
used in the construction of the Novikov complex. 
In the paper \cite{mark}, Mark introduced 
a new class of gradient flows for circle-valued Morse maps 
({\it symmetric flow\/}), 
which are not transversal but, 
somewhat unexpected, 
the Morse-Novikov theory 
can be extended to this case. 
He used these flows to give a yet another 
proof of the Meng-Taubes theorem 
(see the original paper of Meng and Taubes \cite{meng} 
and the later works of Turaev \cite{turaev} and Donaldson \cite{donaldson} 
for alternative proofs of the theorem). 

The symmetric flows have a simple geometric structure allowing 
to carry over to this setting  a large part of the Morse-Novikov theory, 
and on the other hand to perform explicit computations 
with these flows. 
This is the main aim of the present paper. 
We begin by studying the geometric properties of 
symmetric gradients 
(we work actually with a slightly wider class of vector fields 
 called {\it half-transversal gradients\/}), 
and establish the basic theorem of the Morse-Novikov 
theory for this class of flows. 
This theorem is valid in a more general context than 
the Mark's results, 
and we believe that our proof is simpler. 

Then we proceed to detailed study of the geometry 
of the Morse map $f$. 
In the case when $f$ is a fibration 
the first return map from a regular fiber to itself is 
a diffeomorphism, called {\it the monodromy of the fibration\/}; 
this is the basic notion which helps to understand the dynamics of the gradient flow. 
We generalize this notion to the case when $f$ has critical points. 
Our monodromy is a diffeomorphism of two surfaces 
constructed from a Heegaard splitting for the complement 
of a knot \cite{godaone} 
(we recall the basic notions of the theory of 
Heegaard splittings in Section \ref{s:hega}). 
This diffeomorphism depends on the choice of the gradient, 
however it can be efficiently computed 
in particular cases, which leads to the computation 
of the Lefschetz zeta function of certain symmetric gradients 
for the twist knots and the pretzel knot of type $(5,5,5)$. 
The monodromy enables us also to compute the determinant 
of the boundary operator in the Novikov complex for the 
case of these knots 
(the so-called {\it Novikov torsion\/}). 

The dynamics of the gradient flows of 
circle-valued Morse maps are closely related to the 
Seiberg-Witten invariants of 3-manifolds. 
Meng and Taubes \cite{meng}
showed that the Seiberg-Witten invariant of any closed 3-manifold $M$ 
with $b_1(M)\ge 1$ can be identified with the Milnor torsion. 
Fintushel and Stern \cite{f-s} proved that 
for any knot $K$ in $S^{3}$ the Seiberg-Witten invariant 
of the manifold $M\times S^{1}$, where 
$M$ is the result of the zero-surgery on $K$, 
equals the Alexander polynomial of $K$ multiplied 
by a certain standard factor. 
In \cite{donaldson}, Donaldson gives a new proof 
of the Meng-Taubes theorem 
by applying the ideas from Topological Quantum Field Theory. 
These TQFTs were used by Mark to prove a conjecture 
of Hutchings-Lee concerning the relation of the Seiberg-Witten 
invariants with the Novikov torsion.
Some results in this paper have been announced in \cite{godapajitnov2}.

%--------------Section2.tex-----'07/09/18------

\section{Half-transversal flows}
\label{s:half_tr}

Let $f:M\to S^1$ be a Morse function on a closed manifold
$M$.
The dynamics of the gradient flow of $f$ is best understood when
$f$ does not have critical points. In this case 
we choose a regular surface  for $f$,
and the dynamics of the gradient flow is determined by the first return map of 
this surface to itself. 
This map is called   {\it the monodromy}
of the gradient flow.
If $f$ has critical points the situation is much more complicated since
for every transversal $f$-gradient the first return map
is not everywhere defined. It turns out however that in the case of 3-dimensional manifolds
there is an important class of non-transversal gradient flows for which the first return map
determines a self-diffeomorphism of the level surface.
We will first give a definition of the corresponding class of gradient flows on cobordisms.

Let $Y$ be a 3-dimensional cobordism; denote $\pr_-Y,\  \pr_+Y $
the lower, respectively the upper 
boundary of $Y$. Let $\psi:Y\to [a,b]$ 
be a Morse map without critical points of indices 0 and 3.
The subset $U_1$ of all points $x$ in the upper boundary $\pr_+Y $
such that the $(-v)$-trajectory starting at $x$ reaches 
the lower boundary $\pr_-Y $ is  open in $\pr_+Y$
and the gradient descent determines a diffeomorphism 
$(-v)^{\rightsquigarrow}:U_{1}\xrightarrow{\approx}U_{0}$ 
of $U_{1}$ onto an open subset $U_{0}\subset \pr_-Y$.

For two critical points $p,q$ of $f$ 
we call a {\it flow line of $v$ from $q$ to $p$\/}
an integral curve $\gamma$ of $v$ such that 
$$\lim_{t\to-\infty}\gamma(t)=q,\hspace{1cm}
  \lim_{t\to\infty}\gamma(t)=p.$$ 
We shall identify two flow lines of $v$ which are obtained 
from each other by a reparameterization.

\begin{definition}
 A $\psi$-gradient $v$ is called a 
{\it smooth descent gradient\/} 
if
\begin{enumerate}
\renewcommand{\labelenumi}{(\roman{enumi})}
\item
the number of critical points of index 1 is equal to the 
number of critical points of index 2, and they can be 
arranged in two sequences
$$S_1(\psi)=\{p_1,\ldots ,p_k\},\hspace{0.5cm}
  S_2(\psi)=\{q_1,\ldots, q_k\}$$
in such a way that for every $i$ there are two flow 
lines of $v$ joining $q_i$ with $p_i$ and these $2k$ 
flow lines are the only flow lines of $v$. 
\item
the map $(-v)^{\rightsquigarrow}: U_1\to U_0$
can be extended to a $C^{\infty}$ map 
$g:\pr_+Y\to\pr_-Y$.
\footnote{ It seems to us that the point i) actually follows from ii), but we can not prove it at present. }
\end{enumerate}
\end{definition}
Now let us return to circle-valued Morse maps.
Let $f:M\to S^1$ be such a map,
where $M$ is a 3-dimensional closed manifold and $v$ be an $f$-gradient.
Cutting $M$ along a regular surface $S$ of $f$ we obtain a 
cobordism $Y$, a 
Morse function $\psi:Y\to [0,1]$ and a $\psi$-gradient $\bar v=v|Y$.
\begin{definition}\label{d:half_trans}
 The $f$-gradient $v$ is called {\it half-transversal} if there is a regular
level surface $S$ such that $\bar v=v~|~Y$ 
is a smooth descent
 gradient 
of $\psi=f~|~Y$ 
and we have the following transversality condition for stable and unstable manifolds:
\begin{equation}
\mathcal W^{st}(q)\pitchfork \mathcal W^{un}(p)\hspace{0.5cm} 
\end{equation}
for every critical points $p,q$ of $f$ 
with $\ind q=2, \ind p =1$. 
\end{definition}

It is not difficult to show that the subset of all half-transversal gradients 
is dense in the set of smooth descent gradients.

\begin{definition}\label{def:monodromy}
Let $v$ be a half-transversal gradient for a Morse function $f:M\to S^1$
and $S$ be the corresponding level surface of $f$.
The first return map for $(-v)$ determines a diffeomorphism of $S$
to itself which will be  called {\it the monodromy} of the flow generated by $v$, and  denoted by $g$.
\end{definition}

The notion of half-transversal gradient, introduced above originates 
from the paper of T. Mark \cite{mark}
where the class of {\it symmetric flows} was introduced. 
In our terminology Mark's symmetric gradient on a cobordism $Y$
is a smooth descent gradient with the following additional restriction:
there is an involution $I:Y\to Y$ swapping the lower and upper 
boundaries of $Y$ and such that
$I_*(v)=-v$
and $\psi\circ I$ equals $-\psi$ up to an additive constant.
We do not know if the class of smooth descent gradients is 
really wider than Mark's class of symmetric gradients.
However the existence of the involution 
$I$ seems restrictive and we prefer to work with more general notion 
of smooth descent gradients.

Now we will define the Novikov complex and the Lefschetz zeta
function for half-transversal gradient flows. 
The  usual procedure of counting flow lines 
yields the Novikov incidence coefficient 
$$N(q_{i},p_{j};v)=\sum_{k\in\NNN}n_{k}(q_{i},p_{j};v)t^{k}
 \hspace{0.5cm} \in \ZZZ[[t]]$$
where 
$$n_{k}(q_{i},p_{j};v)=\sum_{\gamma\in\Gamma_{k}(q_i,p_{j};v)}\varepsilon(\gamma)$$
(here $\G_k(q_i, p_j; v)$ stands for the set of all flow line of $(-v)$ joining $q_i$ with $p_j$
and $\varepsilon(\gamma)$ is the sign attributed to each flow 
line with respect to the choice of orientations of the 2-dimensional 
stable manifolds). 
The Novikov incidence coefficients
form a square  matrix $D$ with entries 
in $\ZZZ[[t]]$. The chain complex 
\begin{equation}\label{eqn:chain}
0\longleftarrow \mathcal N^-_1\overset{D}{\longleftarrow}\mathcal N^-_2\longleftarrow 0
\end{equation}
where 
$\mathcal N^-_{i}$ is the free $\ZZZ[[t]]$-module 
generated by critical points of $f$ of index $i$ is called 
the {\it positive Novikov complex\/} of the pair $(f,v)$ 
and denoted by $\mathcal N^-_{*}(f,v)$
or simply $\NN^-_*$ if no confusion is possible.
The chain complex
\begin{equation}\label{eqn:chain}
0\longleftarrow \mathcal N_1\overset{D}{\longleftarrow}\mathcal N_2\longleftarrow 0
\end{equation}
where 
$\mathcal N_{i}$ is the free $\ZZZ((t))$-module 
generated by critical points of $f$ of index $i$ is called 
the {\it Novikov complex\/} of the pair $(f,v)$ 
and denoted by $\mathcal N_{*}(f,v)$
or simply $\NN_*$ if no confusion is possible.
The first of the two chain complexes above is more convenient in computations,
however only the homotopy type of the second one is a homotopy 
invariant of the map $f: M\to S^1$ 
(see Theorem \ref{thm:homotopy}).

\begin{definition}
The element $\det D\in \ZZZ[[t]]$
is called {\it the Novikov torsion}
of the pair $(f,v)$ and denoted by 
$\tau(f,v)$ or $\tau_{g}$.
\end{definition}

Proceeding to the Lefschetz zeta functions,
we will need to impose one more
restriction on the gradient flow.
\begin{definition}
Let $f:M\to S^1$ be a Morse function
on a  closed manifold $M$ 
and $v$ an $f$-gradient. 
We say that $v$ is {\it of finite dynamics}
if for every $n\in\ZZZ$ the set of all
closed orbits $\g$
satisfying  $f_*([\g])=n\in H_*(S^1)$
(where $[\g]\in H_1(M)$ is the homology class of $\g$)
is finite. 
\end{definition}

For a half-transversal $f$-gradient of finite dynamics 
we can define the {\it dynamical Lefschetz zeta function of $(-v)$}:
$$\zeta_{-v}(t)=\exp\Big(\sum_{\gamma}
\frac{\varepsilon(\gamma)}{m(\gamma)}t^{m(\gamma)}\Big)$$
where the sum is extended over the set of all closed orbits $\gamma $
of $(-v)$, $\varepsilon(\gamma)$ 
is the Poincar\'e index of $\gamma$, and 
$m(\gamma)$ is the multiplicity of $\gamma$. 
It is clear that  $\zeta_{-v}$ is  equal to the Lefschetz 
zeta function of the diffeomorphism $g$: 
\begin{equation}\label{equ:monodromy}
\zeta_{g}(t)=\exp\Big(\sum_{n\ge 1}\frac{L(g^n)}{n}t^{n}\Big)
\end{equation}
where $L(g^n)$ is the graded trace of the homomorphism induced 
by $g$ in the homology.

Let us now define the class of gradient flows with which we will be working in this paper.
\begin{definition} 
 \label{d:regular}
Let $M$ be a three-dimensional closed manifold, 
and  $f:M\to S^1$ a Morse function 
without critical points of indices 0 or 3;
let $v$ be a half-transversal $f$-gradient of finite dynamics.
We say that $(f,v)$
is a {\it regular  Morse pair}.
\end{definition}

We will also work with Morse functions $f:M\to S^1$ 
on manifolds with
boundary. The definition of the regular Morse pair $(f,v)$
is carried over to this setting in an obvious way,
with the following modifications:

\been\item
The restriction $f~|~\pr M:\pr M\to S^1$ is required to be
a fibration whose  monodromy is isotopic to identity.
\item The gradient vector field $v$ 
is required to be tangent to $\pr M$.
Such gradient  is called {\it a gradient of finite dynamics}
if for every $n\in\ZZZ$ the set of all
closed orbits $\g$  satisfying  $f_*([\g])=n$
is finite. 
\enen

For a regular Morse pair $(f,v)$ on a 3-dimensional manifold with boundary
we define the Novikov complex $\NN_*(f,v)$ and the Lefschetz zeta function
$\zeta_{-v}\in \ZZZ[[t]]$, which counts the closed orbits
of $(-v)$
not belonging to the 
boundary $\pr M$.

%------------------Section3.tex-----07/09/18--------------

\section{The Novikov complex and the 
zeta function of  half-transversal flows}
\label{s:n_z}

The attractive feature of half-transversal flows 
is that the Novikov boundary operators and the Lefschetz zeta function 
of the gradient flow are accessible here through calculations 
with homotopical quantities associated with the monodromy.
Let $M$ be a closed 3-manifold and $(f,v)$  a regular Morse pair 
on $M$.
Let $\overline{M}$ denote the infinite cyclic covering 
of $M$ corresponding to $f$ and $\Delta_{*}(\overline{M})$ 
denote the simplicial chain complex of $\overline{M}$.
Set $\Lambda=\ZZZ[t,t^{-1}]$ and 
$\widehat{\Lambda}=\ZZZ[[t]][t^{-1}]=\ZZZ((t))$.
Both $\NN_*(f,v)$ and $\wh\Delta_{*}(\overline{M})=  
\Delta_{*}(\overline{M})\tens{\L}\wh\L$
are based free finitely generated chain complexes over $\wh\L$.
The next theorem asserts in particular that there is a  chain equivalence between them.
A usual procedure allows to associate to each such equivalence its {\it torsion}, which is
an element in
$$
\Wh(\wh\L)
=
K_1(\wh\L)/U
$$
where $U$ is the subgroup of all elements of the form
$\pm t^n$.
The  group $\Wh(\wh\L)$
is easily identified with the multiplicative group of all
power series in $\ZZZ[[t]]$
with first coefficient equal to $1$
(see \cite{pajitnov06}  Chapter 13, \S 4 for details), so we shall consider the torsions as 
power series with coefficients in $\ZZZ[[t]]$.
The next theorem is the main aim of this section.

\begin{theorem}\label{thm:homotopy}
Let  $M$ be a closed 3-manifold and $(f,v)$  a regular Morse pair 
on $M$.
There is a chain homotopy equivalence 
$$\phi:\mathcal N_{*}(f,v)\to\Delta_{*}(\overline{M})\underset{\Lambda}{\otimes}\widehat{\Lambda}$$
such that 
$$\tau(\phi)=\zeta_{-v}.$$
\end{theorem}

Observe that this theorem implies the isomorphism
$$H_{*}(\mathcal N_{*}(f,v))\approx H_{*}(\overline{M})\underset{\Lambda}{\otimes}\widehat{\Lambda}.$$
Let us first outline the proof.
Lift  $f:M\to S^1$ to a Morse function
$F:\ove M\to\RRR$. The regular level surface $S\sbs M$
(see Definition \ref{d:half_trans}) 
lifts to a regular level surface of $F$ which will be denoted by the same letter $S$.
Denote by $S^{-}$ the part of $\overline{M}$ 
lying below $S$ with respect to 
the function $F$. 
We will construct a certain chain complex $\mathcal Z_{*}$ 
which is free over $\ZZZ[t]$ and computes 
the homology of $S^{-}$. 
Then we construct an embedding 
$$\mathcal N_{*}(f,v)\hookrightarrow
\widehat{\mathcal Z}_{*}=\mathcal Z_{*}
\underset{P}{\otimes}\widehat{P}, 
\hspace{0.5cm}\text{where}\hspace{0.3cm}
P=\ZZZ[t],\,\widehat{P}=\ZZZ[[t]],$$
such that the quotient complex is acyclic and its 
torsion is equal to the Lefschetz zeta function of $-v$. 
The schema of the argument resembles that of the papers 
\cite{hulee2} and \cite{pajitnov99}, however 
the present case is in a sense simpler, 
due to a very particular nature of the half-transversal flows. 

\begin{figure}
\centering
\includegraphics[width=.5\textwidth]{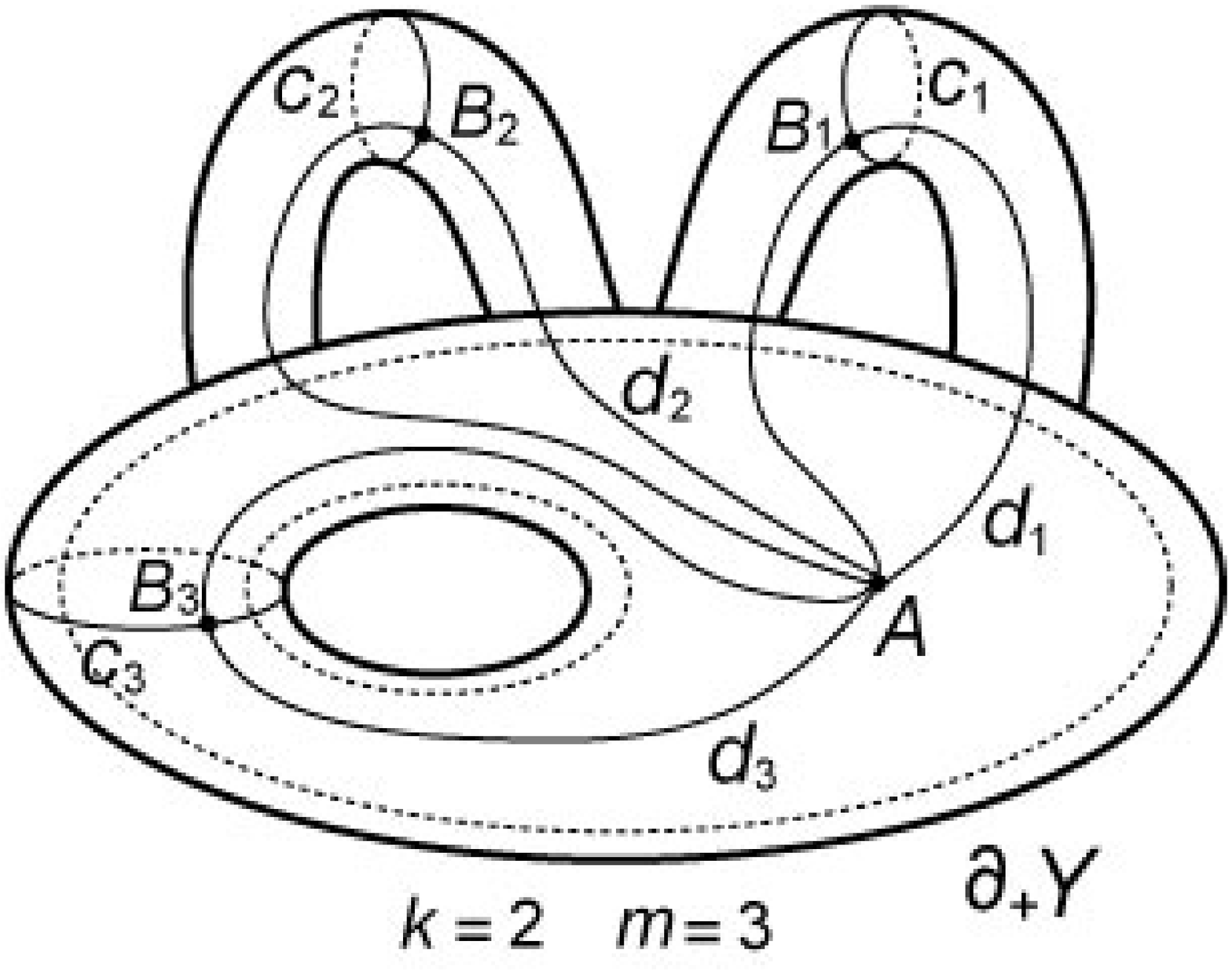}
\caption{}\label{fig:compression}
\end{figure}

Proceeding to details, let us first return to the 
cobordism $Y$ obtained from $M$ by cutting along  $S$.
We have  naturally arising diffeomorphisms $\psi_+:\pr_+Y\to S,\ \ \psi_-:\pr_-Y\to S$.
Put 
$$c_{i}=\WW^{un}(p_{i},v)\cap\pr_+Y. $$
Replacing $Y$ by a diffeomorphic cobordism 
if necessary, 
we can always assume that the circles $c_{i},$ 
$1\le i \le k$ are standardly embedded in $\pr_+Y$ 
as shown in Figure \ref{fig:compression}. 
They are therefore a part of the standard cellular decomposition 
of $\pr_+Y$ which consists of $m$ disjoint circles $c_{i}$, 
and $m$ circles $d_{i}$ 
having a common point $A$. 
For a subset $X\subset\pr_+Y$ we denote $TX$ 
the track of $X$, that is, 
$$TX=\{\gamma(x,t;-v)~|~t\ge 0 \text{ and } x\in X\}.$$

\begin{figure}
\centering
\includegraphics[width=.8\textwidth]{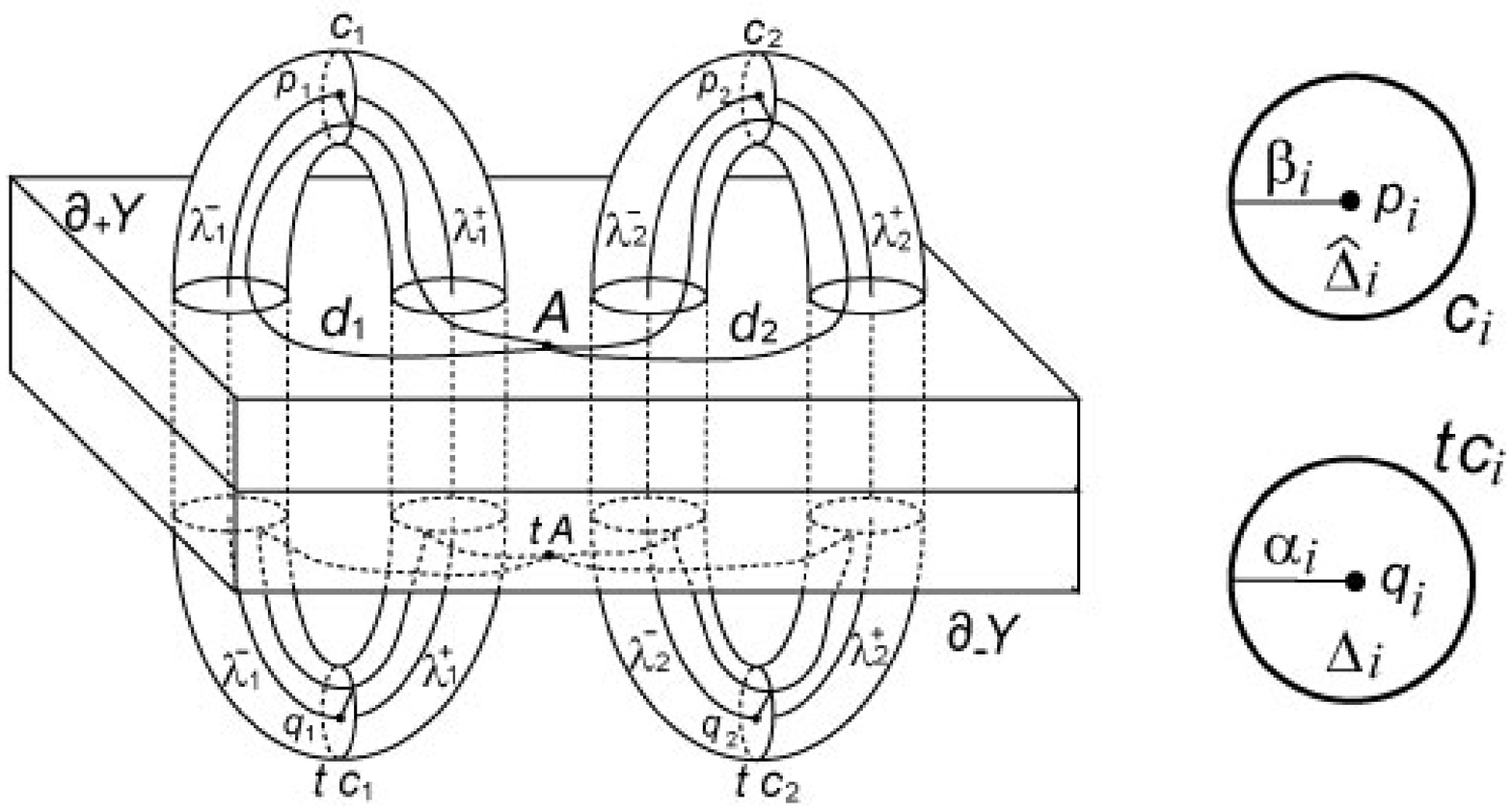}
\caption{}\label{fig:symmetricflow}
\end{figure}

We will now define a filtration $\mathcal E^{i}$
in the cobordism $Y$. 
The term $\mathcal E^{0}$ of the filtration contains two points: $A$ and $tA$. 
The term $\mathcal E^{1}$ contains $\mathcal E^{0}$ 
and the following subsets: 
the circles $d_{i},\,c_{i},$ the 
track $TA$ of the point $A$, the circles 
$$\gamma_{i}=\mathcal W^{un}(q_{i})\cup \mathcal W^{st}(p_i)\hspace{0.5cm}
\text{ for }\hspace{0.3cm}
i\le i\le k$$
the arcs $\aaa_{i}, \,\bbb_{i}$ as shown in Figure \ref{fig:symmetricflow}, 
the circles $Ic_{i},\,Id_{i}\subset \pr_-Y$. 
The term $\mathcal E^{2}$ contains $\mathcal E^{1}$ 
and the following subsets: 
the boundary $\bdd Y$ of $Y$, the stable manifolds of the critical points 
of index 2 and the unstable manifolds of the critical points 
of index 1, and 
the closure of the tracks of $c_{i}$ and $d_{i}$. 
The term $\mathcal E^{3}$ is the whole $Y$.
It is not difficult to see that $\mathcal E^{i}$ is a cellular filtration 
of $Y$, that is, the homology of the quotient $\mathcal E^{i}/\mathcal E^{i-1}$ 
does not vanish only in degree $i$. 

Now we shall use this filtration to explore the homotopy 
type of the covering $\overline{M}$. 
The natural map $Y\to M$ lifts to an embedding of $Y$ to $\ove M$ 
whose image will be identified with $Y$. The covering $\ove M$
is the union of the images of $Y$ under the action of $\ZZZ$:
$$\overline{M}=\underset{n\in\ZZZ}{\bigcup} t^n Y$$
where 
$t$ is the downward generator of $\ZZZ$,
so that $F(tx)=F(x)-1$ for every $x\in \ove M$.
The neighbor copies 
$t^nY$ and $t^{n+1}Y$ are intersecting by 
$\pr_-t^nY=t^n\pr_-Y=t^{n+1}\pr_+Y =\pr_+t^{n+1}Y$.
Recall from Section \ref{s:half_tr}
that the gradient descent determines a diffeomorphism
$g:\pr_+Y\to \pr_-Y$.
We endow $\pr_-Y$ with the cellular decomposition induced from $\pr_+Y$ 
by $g$.
Let $h$ be any cellular approximation of the map 
$\psi_+\circ\psi_-^{-1}:\pr_-Y\to\pr_+Y$.
Then $\ove M$ has the homotopy type of the space
$$
N=\Big(\bigsqcup_{n\in\ZZZ}  t^nY\Big)\Big/\RR
$$
where the equivalence relation $\RR$ 
identifies 
$\pr_-t^nY\approx \pr_-Y$
with
$\pr_+t^{n+1}Y\approx \pr_+Y$
via the map
$h:\pr_-Y\to\pr_+ Y$.
The space $N$ has a natural free action of $\ZZZ$ and we have
a homotopy equivalence
$\ove M\to N$ respecting this action.
Put 
$$
N^-=\Big(\bigsqcup_{n\in \NNN} t^nY \Big) \Big/ \RR.
$$
We will now use the filtration $\mathcal E$ of $Y$ 
to construct a filtration of $N^{-}$. 
Put 
$$\mathcal F^{i}=\underset{n\in\NNN}{\bigcup}t^n\mathcal E^{i}.$$
The filtration $\mathcal S_*(\mathcal F^{i})$ of the singular chain complex $\mathcal S_*(N^-)$ of 
$N^{-}$ is cellular and the homology
$$H_{i}(\mathcal F^{i}/\mathcal F^{i-1})$$
is a free $P$-module. 
Now we will describe the generators of this module.
We denote the stable manifold of $p_{i}$  by $D(p_{i};v)$. 
The set $D(p_{i};v)\setminus\{p_{i}\}$ consists of two arcs, 
their closures will be denoted by 
$\lambda^{+}_{i},\,\lambda^{-}_{i}$ 
(the signs correspond to the chosen orientations). Put
$\lambda_{i}=\lambda^{+}_{i}\cup\lambda^{-}_{i}$.
Let $\b_i$ be an arc in $\D_i$ joining $p_i$ and $B_i=c_i\cap d_i$.
Similarly let $\a_i$ be an arc joining $tA$ with $tB_i$.
Let $d'_{i}$ be the part of $d_{i}$ between  $A$ and $B_{i}$ 
and denote by $\chi^{+}_{i}$ 
the following composition of arcs 
$$\chi^{+}_{i}=d'_{i}\cdot\bbb_{i}\cdot\lambda^{+}_{i}\cdot\aaa_{i}\cdot(td'_{i})^{-1}
\hspace{0.3cm}\text{ where } 1\le i\le k.$$
Similarly, set 
$$\chi^{-}_{i}=d'_{i}\cdot\bbb_{i}\cdot\lambda^{-}_{i}\cdot\aaa_{i}\cdot(td'_{i})^{-1}
\hspace{0.3cm}\text{ where } 1\le i\le k.$$
The fundamental class of $\pr_+Y$ modulo the union of 
$c_{i}$ and $d_{i}$ is denoted by $\omega_{2}$. 
The fundamental class of $Y$ modulo the subspace $\mathcal E^{2}$ 
is denoted by $\omega_3$. 
Here is the list of the free generators of 
$\mathcal Z_{r}=H_{r}(\mathcal F^{r}/\mathcal F^{r-1}):$
as a $\ZZZ[t]$-module:
\begin{alignat*}{2}
r=0 :  & \hspace{0.3cm} A   \\
r=1 :  & \hspace{0.3cm} c_{i},\,d_{i} \hspace{0.5cm} & \text{ for } &1\le i\le m=\mx{ genus} (\pr_+ Y),\\
       & \hspace{0.3cm} \chi^{+}_{i},\,\chi^{-}_{i} \hspace{0.5cm} & \text{ for } &1\le i\le k. \\
r=2 :  & \hspace{0.3cm} \omega_{2}, & \\
       & \hspace{0.3cm} \widehat{\Delta}_{i},\,\Delta_{i},\,Td_{i} \hspace{0.5cm} & \text{ for } &1\le i\le k, \text{ and } \\
       & \hspace{0.3cm} Tc_{i},\,Td_{i}  & \text{ for } &k+1\le i\le m. \\
r=3 :  & \hspace{0.3cm} \omega_3=T\omega_{2}.
\end{alignat*}
Here 
$\widehat{\Delta}_{i}$ 
is the unstable manifold of $p_i$ in $Y$;
we have $\pr\widehat{\Delta}_{i}=c_i$,
and similarly for $\D_i$.
(By a certain abuse of notations we use the same symbol $c_{i}$ for the cycle 
 and its geometric support; 
 similar convention holds for the other notations.)
Now we shall describe 
the boundary operators in the adjoining complex 
$$\bdd_{r} : \mathcal Z_{r}\to \mathcal Z_{r-1}:$$ 

$\bdd_{1} : \mathcal Z_{1}\to\mathcal Z_{0} :$
$$\bdd(c_{i})=0=\bdd(d_{i}),\hspace{0.3cm} \bdd(\chi^{+}_{i})=\bdd(\chi^{-}_{i})=\bdd(TA)=A-th(A).$$

$\bdd_{2} : \mathcal Z_{2}\to\mathcal Z_{1} :$
\begin{equation*}
\left.
\begin{aligned}
\bdd(\widehat{\Delta}_{i}) & = -c_{i} \\
\bdd(\Delta_{i}) & = th(c_{i}) \\
\bdd(Td_{i}) & = d_{i}+\lambda_{i}-th(d_{i}) 
\end{aligned}
\hspace{0.5cm}\right\} \text{ for } 1\le i\le k, 
\end{equation*}

\begin{equation*}
\left.
\begin{aligned}
\bdd(Tc_{i}) & = c_{i}-th(c_{i}) \\
\bdd(Td_{i}) & = d_{i}-th(d_{i})
\end{aligned}
\hspace{0.5cm}\right\} \text{ for } k+1\le i\le m, \text{ and }
\end{equation*}

$$\bdd(\omega_{2})=0.$$

$\bdd_{3} : \mathcal Z_{3}\to\mathcal Z_{2} :$
$$\bdd(\omega_3)=\omega_{2}-th(\omega_{2}).$$
The chain complex $\mathcal Z_{*}$ is chain equivalent to the 
simplicial chain complex of $N^{-}$.
Any chain equivalence 
$$
\xi:\ZZ_*\to\D_*(N^-)
$$
has a well-defined torsion $\tau(\xi)\in \Wh(\ZZZ[t])
=
K_1(\ZZZ[t])/\{\pm 1\}$.
This last group vanishes
(by the Bass-Heller-Swan theorem), therefore $\tau(\xi)=0$, and
the torsion of the chain equivalence 
$$
\wh\xi:\wh\ZZ_*=\ZZ_*\tens{\ZZZ}\ZZZ[[t]] \to \D_*(N^-)\tens{\ZZZ}\ZZZ[[t]]
$$
in the group 
$K_1(\ZZZ[[t]])/\{\pm 1\}$
vanishes. To prove our theorem it suffices 
therefore to construct a chain equivalence
$$
\NN_*^-=\NN_*^-(f,v)\overset{\sigma}\rightarrow \wh\ZZ_*
$$
such that
$\tau(\sigma)=\zeta_{-v}$. We will embed 
$\NN_*^-$ to $\wh\ZZ_*=\ZZ_*\tens{\ZZZ[t]}|\ZZZ[[t]]$
and compute its quotient complex.
Let us first observe that the Novikov complex 
for our half-transversal flow can be
expressed in terms of the monodromy $g$ or its
homotopy substitute $h$:
$$
\pr q_i = \sum N(q_i, p_j) p_j,
\quad
\mx{ where }
N(q_i, p_j)
=
\sum_{k\in\NNN} t^k\langle h^k(c_i), c_j\rangle
$$
where $\langle\cdot, \cdot \rangle $
stands for the pairing in $H_1(\pr_+Y)$.
Now let us make a simple change of basis 
\footnote{ A change of basis is called {\it simple}
if the torsion of the transition matrix 
vanishes in $\Wh\big(\ZZZ[[t]]\big)
=
K_1(\ZZZ[[t]])/\{\pm 1\}$.} 
in 
$\mathcal Z_{*}$ replacing $\widehat{\Delta}_{i}$ 
by the element $\widehat{\Delta}_{i}-\Delta$ 
which will be denoted by $Tc_{i}$ (in order to stress the analogy with 
the tracks of the circles $d_{i}$). 
Extending the map $T$ by linearity to a homomorphism 
$H_{1}(\pr_+Y)\rightarrow \mathcal Z_{2}$ 
it is easy to check the following formula: 
\begin{equation}\label{equ:base1}
\bdd(T\mu)=\mu-th(\mu)+\sum_{j}\langle\mu,c_{j}\rangle\lambda_{j}.
\end{equation}

Let us now make one more simple change of basis, 
replacing the cycle $\Delta_{i}$ by 
\begin{equation}\label{equ:base2}
\widetilde{\Delta}_{i}=\Delta_{i}-\sum_{j=1}^{\infty}t^{j}T(h^{j}c_{i}).
\end{equation}

This infinite sum corresponds geometrically to the stable manifold 
of the critical point $p_{i}$. 
There is however one essential difference between the formula (\ref{equ:base2}) 
and the similar formulas for the case of the transversal flows 
(see, for example, formula (66) from \cite{pajitnov99}). 
The formula (\ref{equ:base2}) contains the term 
$Tc_{i}=\widehat{\Delta}_{i}-\Delta_{i}$ 
and similar ones
which are not strictly speaking the geometric traces of the cells. 
An easy computation using the formula (\ref{equ:base1}) 
shows that the 
homomorphism 
$\s:\NN_*^-\to\ZZ_*$ defined by
$$
\s(p_i)=\l_i, \ \ \s(q_i)=\wi\D_i
$$
is an embedding of chain complexes.
The quotient complex $Q_{*}$ is also 
easily computed; 
here is the list of free $\ZZZ[[t]]$-generators for $Q_{j}$: 
\begin{align*}
j=0 :& \hspace{0.7cm} A   \\
j=1 :& \hspace{0.7cm} TA,\,c_{i},\,d_{i},\,\chi^{+}_{i}\\
j=2 :& \hspace{0.7cm} Tc_{i},\,Td_{i},\,Td'_{i},\,\omega_{2}\\
j=3 :& \hspace{0.7cm} T\omega_{3}
\end{align*}
We have $\bdd(Td'_{i})=\chi^{+}_{i}$ 
and 
$$\bdd(z)=0, \hspace{0.3cm}\hspace{0.3cm} \bdd(Tz)=1-th(z)$$ 
for every $z$ from the following list:  
$$A,\,c_{i},\,d_{i},\,\omega_{2}.$$ 
 
After factoring out the chain complex generated by $\chi^{+}_{i}$ 
and $Td'_{i}$, we obtain the chain complex of the mapping torus 
of the map $h$. 
It is well known that its torsion equals the Lefschetz
zeta function of $h$ 
(see the classical paper of J. Milnor \cite{milnor2}). 
This completes the proof of Theorem \ref{thm:homotopy}.
$\qs$
	
\begin{remark}\label{r:nov_zeta_boundary}
The theorem above is valid also in the case 
of regular Morse pairs on manifolds with boundary, 
and the proof is similar.
\end{remark}

%------------------Section4.tex-----07/09/18--------------

\section{Novikov torsion and the Alexander polynomial for knots}
\label{s:tors_alex}

Theorem \ref{thm:homotopy}
establishes a relation between two natural geometric objects: the homotopy equivalence 
$\phi:\NN_*(f,v)\to \D_*(\overline{M})\tens{\L}\wh \L$
and the Lefschetz zeta function of the  flow generated by $v$.
For computational purposes it is convenient to reformulate it 
in another way.
Let $(f,v)$ be a regular Morse pair on a  
3-manifold $M$ (with or without boundary).
Let $\FFF$ be a field. 
\begin{definition}\label{d:def_acyclic}
We say that $(f,v)$ is {\it $\FFF$-acyclic},
if
$$
H_*(\ove M)\tens{\ZZZ[t,t^{-1}]} \FFF((t))=0.
$$
\end{definition}
Put $\NN_*(f,v;\FFF)=\NN_*(f,v)\tens{\wh\L}\FFF((t)).$
It follows from Theorem \ref{thm:homotopy}
that if $(f,v)$ is $\FFF$-acyclic, then the homology
of the complex $\NN_*(f,v;\FFF)$ also vanishes.
The images of the elements $\tau(f,v), \z_{-v}$
in the ring $\FFF[[t]]$ will be denoted by
$\tau^\FFF, \ \z_{-v}^\FFF$.
The $\FFF$-acyclicity condition implies that the torsion of the chain complex 
$$
\wh\D_*^\FFF(\ove M)=\D_*(\ove M)\tens{\ZZZ[t,t^{-1}]} \FFF((t))
$$
is well defined as an element of
$$
\Wh(\FFF((t)))
\approx
K_1(\FFF((t)))/U,
$$
where $U$ is the subgroup of all elements of the form
$\pm t^n$. 
We will denote this torsion by $\tau_M^\FFF$ 
omitting in the notation the obvious dependence of this element
on the homotopy class of $f$.

\begin{proposition}
 In the assumptions of Theorem \ref{thm:homotopy}
assume moreover that 
$(f,v)$ is $\FFF$-acyclic.
Then
$$
\tau^\FFF\cdot \zeta^\FFF_{-v}=\tau^\FFF_M.
$$
\end{proposition}
\Prf
Tensoring by $\FFF((t))$ the chain equivalence $\phi$
we obtain a chain equivalence
$$
\phi^\FFF:\NN_*(f,v;\FFF)\to 
\wh\D^\FFF(\ove M)
$$
of two acyclic complexes.
The torsion of such chain equivalence equals the quotient of the torsions of the complexes. $\qs$

Let $K\sbs S^3$ be an oriented  knot, $M=S^3\sm \Int N(K)$,
and $\FFF=\QQQ$. Let $(f,v)$ be a regular Morse pair
on $M$ such that the homotopy class
$[f]\in H^1(M)\approx [M, S^1]\approx \ZZZ$
is the positive generator of this group.
The condition of $\QQQ$-acyclicity is fulfilled here,
so the above proposition is valid.
It is well known that in this case the torsion $\tau_M$
equals the Alexander polynomial divided by $(1-t)$
and we obtain the following corollary:

\begin{corollary}
 Let $K$ be a knot in $S^3$, let $M=S^3\sm \Int N(K)$
and $(f,v)$ be a regular Morse pair on $M$.
 Let 
$\tau$ be the Novikov torsion of $(f,v)$.
Then 
$$
\tau\cdot \zeta_{-v}
=
\frac {\D_K}
{1-t}
$$
where $\D_K$ stands for the Alexander 
polynomial of the knot $K$.
\end{corollary}

%--------------Section5.tex-----'07/09/20------

\section{Heegaard splitting for sutured manifolds}
\label{s:hega}

The notion of a sutured manifold was introduced by Gabai \cite{gabai1}. 
See also \cite{s-0}.
In this section, we recall the notations and define Heegaard splitting 
for the sutured manifolds \cite{godaone}. 
%Then we discuss the stabilization problem of this splittings. 

\begin{definition}
A {\it sutured manifold\/} $(X, R_+, R_-)$ is a compact oriented 3-manifold $X$ 
with $\partial X$ decomposed into the union along the boundary of two 
connected surfaces 
$\tilde{R}_+$ and $\tilde{R}_-$ oriented so that 
$\partial\tilde{R}_+=\partial\tilde{R}_-=\gamma$ and 
$\partial X=\tilde{R}_+\cup\tilde{R}_-$.
Let $A(\gamma)$ denote a collection of disjoint annuli comprising 
a regular neighborhood $\gamma$, and define 
$R_{\pm}=\tilde{R}_{\pm}-\text{Int }A(\gamma)$. 
Thus $\partial X=R_+\cup R_-\cup A(\gamma)$. 
We regard $R_+$ as the set of components of $\partial X-\text {Int }A(\gamma)$ 
whose normal vectors point out of $X$, and $R_-$ as those whose normal 
vectors point into $X$.  
The symbol will denote $R_+$ or $R_-$ respectively while $R(\gamma)$ denotes 
$R_+\cup R_-$.
If $\partial\tilde{R}_+=\partial\tilde{R}_-=\emptyset$,   
each component of $\tilde{R}_{\pm}=R_{\pm}$ is a closed surface. 
\end{definition}

Let $L$ be a non-split oriented link in a homology 3-sphere, 
and 
$\bar{R}$ a Seifert surface of $L$. 
Set $R=\bar{R}\,\cap E(L)$ ($E(L)=$ cl$(S^3-N(L))$). 
Let $P$ be a regular neighborhood of $R$ in $E(L)$, then 
$P$ forms $R\times [-1, 1]$ where $R=R\times\{ 0\}.$ 
We denote by $\acute{R}_+$ ($\acute{R}_-$ resp.) $R\times\{ 1\}$ 
($R\times\{ -1\}$ resp.), then 
$(P,\acute{R}_+,\acute{R}_-)$ may be regarded as a sutured manifold. 
We call $(P, \acute{R}_{+}, \acute{R}_{-})$ a {\it product sutured manifold for\/} $R$. 
Further, 
let $X=$ cl$(E(L)-P)$, and $R_{\pm}=\acute{R}_{\mp}$, then 
we may also regard $(X,R_+,R_-)$ as a sutured manifold. 
We call $(X,R_+,R_-)$ the
{\it complementary sutured manifold for\/} $R$. 
In this paper, we call this the {\it sutured manifold for\/} $R$
for short. 

\begin{figure}
\centering
\includegraphics[width=.7\textwidth]{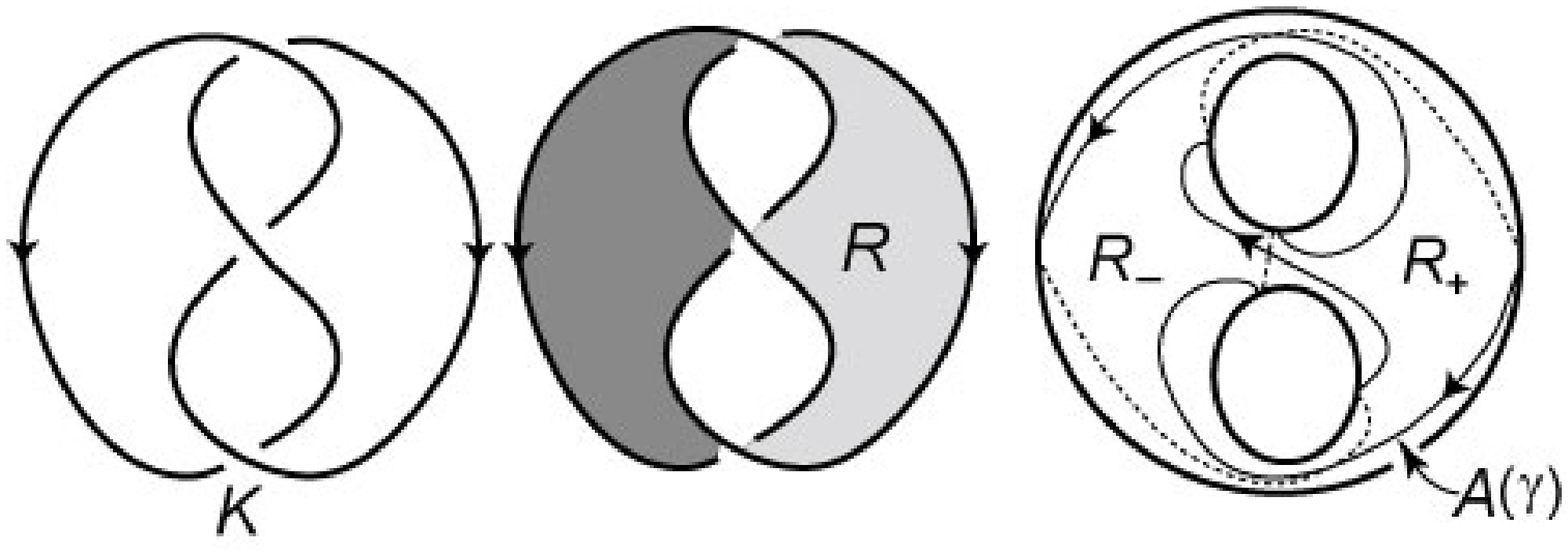}
\caption{}\label{fig:trefoil}
\end{figure}

\begin{example}
Let $K$ be the trefoil knot in the 3-sphere $S^3$ and 
$R$ the genus 1 Seifert surface as illustrated in Figure \ref{fig:trefoil}. 
The (complementary) sutured manifold for $R$ is homeomorphic to 
the manifold in the righthandside of the figure.  
(Note that the `outside' of the genus 2 surface is 
the complementary sutured manifold.) 
\end{example}

\begin{definition} \label{def:compression}
A {\it compression body\/} $W$ is a connected 3-manifold obtained 
from a compact surface $\partial_-W$ by attaching 1-handles to 
$\partial_-W\times \{ 1\}\subset \partial_-W\times [0,1]$. 
Dually, a compression body is obtained from a connected surface 
$\partial_+W$ by attaching 2-handles to $\partial_+W\times\{1\}
\subset\partial_+W\times [0,1]$ and 3-handles to any spheres thereby 
created. 
If $W=\bdd_+W\times [0,1]$, $W$ is called a {\it trivial\/} compression body.   

We collapse a compression body $W$, so that we may obtain 
$\bdd_-W\cup\text{ (arcs)}$, where the arcs correspond to cores 
of the attaching 1-handles. We say the family of arcs 
the {\it spine\/} of $W$. 
We denote by $h(W)$ the number of the attaching 1-handles of $W$.
\end{definition}

\begin{definition}
A pair $(W,W')$ is a Heegaard splitting for a sutured manifold $(X,R_+,R_-)$ if :
\begin{enumerate}
\renewcommand{\labelenumi}{(\roman{enumi})}
\item 
both $W$ and $W'$ are compression bodies; 
\item
$W\cup W'=X$; 
\item
$W\cap W'=\partial_{+}W=\partial_{+}W', 
\partial_{-}W=R_{+}$ and $\partial_{-}W'=R_{-}$. 
\end{enumerate}
\end{definition}

If $\gamma\neq\emptyset$, then $\partial_-W$ and $\partial_-W'$ 
have boundaries so that 
$\bdd(\bdd_-W)\times [0,1] \cup \bdd(\bdd_-W')\times [0,1]=A(\gamma)$ and 
$\bdd(\bdd_+W)=\bdd(\bdd_+W')=\gamma$. 
This case are treated in \cite{godaone} and \cite{godatwo}.  
See also \cite{godathree} for the concrete examples.  
We should note that if $R_+$ is homeomorphic to $R_-$, 
we have $h(W)=h(W')$. 

\begin{remark}\label{rmk:morse}
This Heegaard splitting corresponds to a circle-valued Morse map 
$M\rightarrow S^1$ for a closed orientable 3-manifold 
$M$ 
with $b_1(M)>0$ or the complement of 
a non-split link in a homology 3-sphere $M$. 
In both cases, we suppose that we have a compact surface $R$ 
as a representative of $H_1(M)$. 
Then, we obtain 
the sutured manifold $(X,R_+,R_-)$ from $M$ by cutting along $R$. 
So, we have a Heegaard splitting $(W, W')$ of $(X,R_+,R_-)$ as above. 
%Here, there is a Morse map corresponding to this splitting.
%The index 1 critical points of the Morse map correspond to the attaching 
%1-handles of $W$ and index 2 critical points correspond to those 
%of $W'$. 
See \cite{godapajitnov} and \cite{p-r-w} for the detail.
\end{remark}

\begin{definition}
Suppose that $R_+$ is homeomorphic to $R_-$. 
Set 
$h(X,R_+,R_-)=\min\{h(W)(=h(W'))~|~(W,W') \text{ is a Heegaard splitting for }(X,R_+,R_-)\}.$ 
We call it the {\it handle number \/} of $(X,R_+,R_-)$. 
The {\it Morse-Novikov number \/} $\mathcal M\mathcal N$ of $(M,R)$ or $(X,R_+,R_-)$  
is the minimal possible number of the critical points of 
the corresponding Morse map.
%corresponding to $(W,W')$. 
\end{definition}

\begin{remark}
By Corollary 2.8 in \cite{godapajitnov}, 
we may see that 
$\mathcal M\mathcal N(M,R)=2\times h(X,R_+,R_-).$
\end{remark}

\begin{definition}\label{def:stabilization}
Suppose that $(W,W')$ is a Heegaard splitting of a sutured manifold 
$(X,R_+,R_-)$, and 
let $\lambda$ be a properly embedded arc in $W'$ parallel to 
an arc in $\bdd_+W'$. 
Here ``parallel'' means that there is an embedded disk $D$ 
in $W'$ whose boundary is the union of $\lambda$ and an arc in $\bdd_+W'$. 
Now add a neighborhood of $\lambda$ to $W$ and delete it from $W'$. 
This adds a 1-handle to $W$ (whose core is $\lambda$) and also 
adds a 1-handle to $W'$ (whose cocore is a disk in $D$). 
Thus we have again the Heegaard splitting $(\widehat{W}, \widehat{W}')$ 
of $(X,R_+,R_-)$ where the genus of $\widehat{W}$ ($\widehat{W}'$ resp.) 
is one greater than $W$ ($W'$ resp.). 
This process is called a {\it stabilization\/} of $(W,W')$.
\end{definition}

We may regard a compression body $W$ as a sutured manifold $(W,R_+,R_-)$,
that is, we may suppose $\bdd_+W=R_+$ and $\bdd_-W=R_-$. 
%Assume that $A(\gamma)=\emptyset$.
A compression body $W$ has a natural Heegaard splitting: 
A surface $S$ parallel to $\bdd_+W$ splits $W$ into two compression bodies, 
at least one of them is trivial. 
Call this the {\it trivial splitting\/} of $W$. 
A splitting is called {\it standard\/} 
if it is obtained from the trivial splitting by stabilization.
In \cite{s-t}, Scharlemann and Thompson proved the next theorem: 

\begin{theorem}[\cite{s-t}]\label{thm:standard}
Every Heegaard splitting of a compression body $(W,R_{+},R_{-})$ 
with $\gamma=\emptyset$ is standard.
\end{theorem}

\begin{remark}
In \cite{s-t}, two types of trivial splittings, called `type 1 and 2', are treated. 
Here we have only to consider the `type 1' trivial splitting.   
\end{remark}

This theorem induces the following theorem. 
The idea is due to Lei \cite{lei}. 

\begin{theorem}\label{thm:stabilization}
Any two Heegaard splittings of the the same sutured manifold with $\gamma=\emptyset$ 
have a common stabilization.
\end{theorem}

\begin{proof}
Let $(W,W')$ and $(V,V')$ be Heegaard splitting of a sutured manifold 
$(X,R_+,R_-)$ with $\gamma=\emptyset$ 
such that $\bdd_-W=R_+$  and $\bdd_-V'=R_-$.  
Let $\lambda_W$ and $\lambda_{V'}$ be the spines of $W$ and $V'$. 
Then, the standard general position argument allows that 
$N(\partial_{-}W\cup\lambda_W)\cap N(\partial_{-}V'\cup\lambda_{V'})=\emptyset$. 
We denote by $X$ the sutured manifold with $R_+=\bdd_+W$ and 
$R_-=\bdd_+V'$, and let $S$ be a Heegaard splitting surface for $X$. 
Then $S$ is also a Heegaard splitting surface for $(X,R_+,R_-)$. 
Moreover, $S$ becomes a Heegaard splitting surface for the compression bodies
$W'=X-\text{ Int }N(\partial_{-}W\cup\lambda_W)$ and 
$V=X-\text{ Int }N(\partial_{-}V'\cup\lambda_{V'})$. 
Hence the Heegaard splitting surface $S$ is 
a stabilization of both $(W,W')$ and $(V,V')$ by Theorem \ref{thm:standard}. 
\end{proof}

As in Remark \ref{rmk:morse}, 
if there is a circle-valued Morse map $f : M\to S^1$, 
we have a Heegaard splitting $(W,W')$ of the 
sutured manifold $(X,R_{+},R_{-})$. 
We also say that $(W,W')$ is a Heegaard splitting of $M$
or $Y$.
Let $\lambda_{W}=\cup_{i}\lambda_{W}^{i}$  
($\lambda_{W'}=\cup_{i}\lambda_{W'}^{i}$ resp.) 
be the set of spines of $W$ ($W'$ resp.). 

\begin{definition}\label{def:symmetric}
A family $(W,W',\lambda_{W},\lambda_{W'})$ is called 
a {\it symmetric Heegaard splitting\/} of $M$ if 
it satisfies the following conditions: 
\begin{enumerate}
\renewcommand{\labelenumi}{(\roman{enumi})}
\item
$(W,W')$ is a Heegaard splitting of $M$; 
\item 
there is one to one correspondence 
between the arcs $\lambda_{W}^{i}$ and $\lambda_{W'}^{i}$ 
$(i=1,\ldots, k)$. 
Further, 
$\partial\lambda_{W}^{i}=\partial\lambda_{W'}^{i}$ 
for each $i$.
\end{enumerate}
\end{definition}

\begin{remark}
For a half-transversal gradient flow, 
we can construct a symmetric Heegaard splitting 
so that $\cup_{i}(\lambda_{W}^{i}\cup\lambda_{W'}^{i})$ 
are the circles of the half-transversal flow. 
Conversely, for every symmetric Heegaard splitting $\mathcal H$, 
there is a homeomorphism $\varphi$ of $Y$
such that $\varphi(\mathcal H)$ 
is obtained from a half-transversal gradient flow.
\end{remark}

%--------------Section6.tex-----'07/09/20------

\section{Counting closed orbits}
\label{s:clo_orb}

In this section, we establish a method to count 
closed orbits using the idea described in the previous sections.

Let $R$ be compact connected manifold, $g: R\rightarrow R$ be a continuous map. 
Assume that $g$ has only finite number of the critical points.
The {\it Lefschetz number\/} is defined as follow:
$$L(g)=\sum_{i=1}^{\ell}\textrm{ind}(x_i),$$
where  ind$(x_i)$ is the index of the fixed point $x_{i}$ 
(see \cite{dold}). 
Let $G_i$ be the endomorphism of the homology group $H_i(R)$ 
induced by $g$. 
Then the Lefschetz fixed point theorem asserts the following: 
\begin{equation}\label{equ:trace}
L(g)=\sum_i(-1)^i{\rm trace}(G_i: H_i(R) \rightarrow H_i(R)).
\end{equation}

Let $K$ be a fibred knot in the 3-sphere $S^3$. 
Then $K$ has a Seifert surface $R$ 
and 
the complement of $K$ is the fiber bundle over $S^1$ 
with fiber $R$. 
Let $(P,\acute{R}_{+},\acute{R}_{-})$ be the product sutured manifold for $R$, 
and $(X,R_{+},R_{-})$ the complementary sutured manifold for $R$.  
Then $(X,R_{+},R_{-})$ has also product sutured manifold structure.

The monodromy $g$ induces the transformation matrix 
$G_{i}: H_i(R)\rightarrow H_i(R).$
We call $G_{1}$ the {\it monodromy matrix\/} of the 
fibred knot $K$. 
Concretely, 
we can have a presentation of $G_{1}$ as follows.
Let $c_1,c_2, \ldots, c_{m},d_1,d_2, \ldots, d_{m}$ be 
symplectic basis of $H_1(R)$, 
where $m$ is the genus of $R$. (See e.g. \cite{morita}.)
We suppose that $c_{i}\cdot d_{i}=1$ here.
Push them off along the normal vector of $R$, 
and put them on $\acute{R}_{+}$ and $\acute{R}_{-}$. 
Then we may see that they are basis 
of $H_1(\acute{R}_{+})$ and $H_1(\acute{R}_{-})$. 
Since $R_{\pm}=\acute{R}_{\mp}$, 
we may denote the basis of $H_1(R_{+})$ ($H_1(R_{-})$ resp.) 
by $c^{+}_1,\ldots, c^{+}_{m}, d^{+}_{1},\ldots,d^{+}_{m}$ 
($c^{-}_1,\ldots,c^{-}_{m},d^{-}_{1},\ldots,d^{-}_{m}$ resp.). 
By using the product structure of $(X,R_{+},R_{-})$, 
we push further 
$c^{-}_1,\ldots,c^{-}_{m},d^{-}_{1},\ldots,d^{-}_{m}$ into 
$R_{+}$, and 
denote their images in $H_{1}(R_{+})$ 
by 
$c'_1,\ldots,c'_{m},d'_{1},\ldots,d'_{m}$.
Then, 
$$
\begin{pmatrix}
	c'_1 \\
	c'_2 \\
	\cdot \\
	\cdot \\
	\cdot \\
	d'_{m} 
	\end{pmatrix} 
=G_1
\begin{pmatrix}
	c_1^+ \\
	c_2^+ \\
	\cdot \\
	\cdot \\
	\cdot \\
	d_{m}^+ 
\end{pmatrix}.$$

%This observation can be obtained in case of links in $S^3$.
We show an example here. 

\begin{figure}
\centering
\includegraphics[width=.5\textwidth]{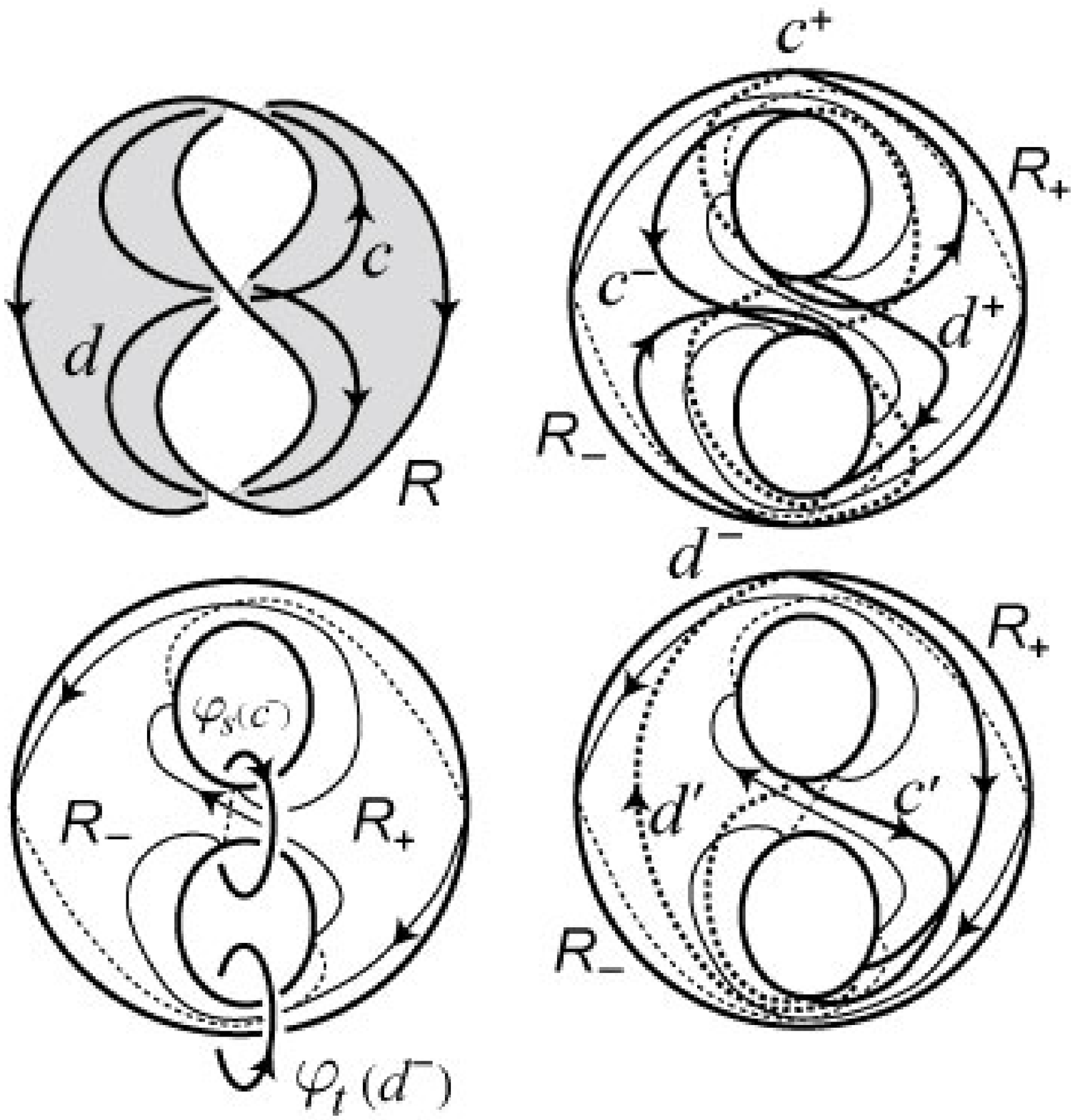}
\caption{}\label{fig:trefoil2}
\end{figure}

\begin{example}\label{ex:trefoil}
Let $K$ be the trefoil knot and $R$ the Seifert 
surface as shown in Figure \ref{fig:trefoil}. 
Set $c$ and $d$ as generators of $R$ illustrated 
in Figure \ref{fig:trefoil2}. 
The upper right-hand figure in Figure \ref{fig:trefoil2} 
shows that the sutured manifold $(X,R_+,R_-)$ for $R$
with $c^{\pm},d^{\pm}\,\subset R_{\pm}$. 
This (complementary) sutured manifold $X$ is a product sutured manifold, 
that is, $X$ is homeomorphic to $R\times [0,1]$ 
where $R_{-}=R\times\{ 0\}$ and $R_{+}=R\times\{ 1\}$. 
Then we can consider a `flow' $\varphi_s$ $(s\in [0,1])$ 
using this product structure
such that $\varphi_{s}(a)=a\times\{ s\}\subset R\times\{ s\}$ 
for a subset $a$ in $R_{-}$.
$\varphi_{s}(c^{-})$ and $\varphi_{t}(d^{-})$
$(s,t \in (0,1), (s\neq t))$ are depicted in the lower left-hand figure in 
Figure \ref{fig:trefoil2}, 
and the lower right-hand figure shows 
$\varphi_1(c^{-})$ and 
$\varphi_1(d^{-})$, 
denoted by $c'$ and $d'$.
Therefore we can observe that 
$$
\begin{pmatrix}
	c' \\
	d' 
\end{pmatrix}
=
\begin{pmatrix}
	d^{+} \\
	-c^{+}+d^{+}
\end{pmatrix}
=
\begin{pmatrix}
	\phantom{-}0 & 1 \\
	          -1 & 1 
\end{pmatrix}
\begin{pmatrix}
	c^{+}  \\
	d^{+}
\end{pmatrix}.
$$
Thus we have 
$$G_1=
\begin{pmatrix}
	\phantom{-}0 & 1 \\
	          -1 & 1 
\end{pmatrix}.$$
\end{example}

In this case, we can observe that $\text{trace}(G_0: H_0(R)\rightarrow H_0(R))=1$ and 
$G_2=0$.
From (\ref{equ:monodromy}) and (\ref{equ:trace}), we have :  
\begin{align*}
\zeta_{g}(t)
&=\exp\big(\sum_{k=1}^{\infty}\frac{t^k}{k}(1-\text{trace } G_1^k)\big) \\
&=\exp\big(\log(1-t)^{-1}+\text{trace}(\log(I-t\cdot G_1))\big) \hspace{0.7cm} (\,|t|<1\,) \\
&=\frac{\det(I-t\cdot G_1)}{1-t} \\
&=\frac{1-t+t^2}{1-t}.
\end{align*}
Here $I$ is the unit matrix. 
Note that the Alexander polynomial of the trefoil knot is 
$1-t+t^2$. In general, if a knot $K$ is fibred, 
the numerator $\det(I-t\cdot G_{1})$ 
equals the Alexander polynomial of $K$.
Therefore we have the following well-known theorem.
See \cite{milnor2} for example. 

\begin{theorem}[\cite{milnor2}]
Let $K$ be a fibred knot in $S^3$, and we denote by
$g$ the monodromy of $K$.  Then, 
$$\zeta_{g}(t)=\frac{\Delta_K(t)}{1-t}.$$
Here $\Delta_K(t)$ is the Alexander polynomial of $K$.
\end{theorem}

Now let us consider the case of non-fibred. 

Let $M$ be a compact orientable 3-manifold 
with $b_{1}(M)>0$. 
Let $f : M\to S^1$ be a Morse map, 
and $R$ a regular level surface for $f$.  
We obtain a sutured manifold $(X,R_+,R_-)$ 
from $M$ cutting along $R$.
As pointed out in Remark \ref{rmk:morse} and 
Definition \ref{def:symmetric}, 
there is a symmetric Heegaard splitting $(W,W',\lambda_{W},\lambda_{W'})$ 
corresponding to $f$. 
Set $k=h(W)(=h(W'))$ the number of the attaching 1-handles 
of $W$. 

\begin{figure}
\centering
\includegraphics[width=.4\textwidth]{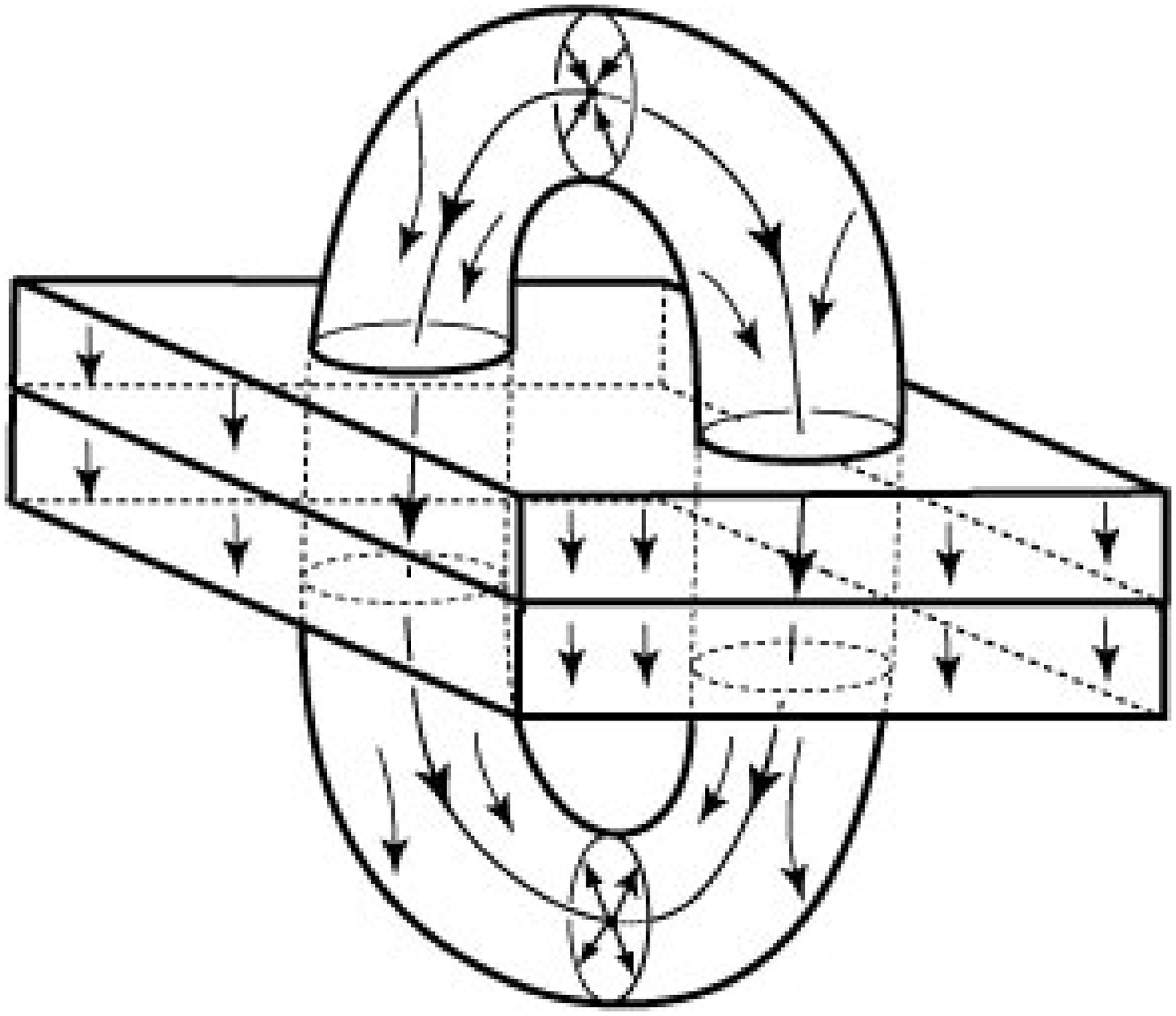}
\caption{}\label{fig:flow}
\end{figure}

According to Definition \ref{def:monodromy}, the monodromy $g$ 
induces the transformation matrix $G_{1}:H_{1}(S)\to
H_{1}(S)$, which can be obtained as follows.
We denote the symplectic basis of $H_{1}(\partial_{+}W)$ 
($H_{1}(\partial_{+}W')$ resp.) by 
$c^{+}_1,,\ldots, c^{+}_{k}, c^{+}_{k+1},\ldots ,c^{+}_{m}$, and 
$d^{+}_1,\ldots,d^{+}_{k},d^{+}_{k+1},\ldots, d^{+}_{m}$ 
($c^{-}_1,\ldots, c^{-}_{k},c^{-}_{k+1},\ldots,c^{-}_{m}$, and  
$d^-_1,\ldots,d^-_{k}, d^-_{k+1},\ldots,d^-_{m}$ resp.). 
Here $c^{+}_j$ and $d^{+}_j$ ($c^{-}_j$ and $d^{-}_j$ resp.) $(j=1,\ldots ,k)$ 
are derived from 
the attaching 1-handles of $W$ ($W'$ resp.), namely, 
$c^{+}_j$ ($c^{-}_j$ resp.) $(j=1,\ldots ,k)$ is a cocore of the attaching 1-handle of 
$W$ ($W'$ resp.) and $d^{+}_j$ ($d^{-}_j$ resp.) $(j=1,\ldots ,k)$ is a `longitude' corresponding to $c^{+}_j$ 
($c^{-}_{j}$ resp.), 
so that $c^{+}_j\cdot d^{+}_{\ell}=\delta_{j\ell}=c^{-}_{j}\cdot d^{-}_{j}$ $(j,\ell=1,\ldots ,m)$.
cf. Figure \ref{fig:compression}. 
As in the case of a fibred knot, 
the generators $c^{+}_{j}, d^{+}_{j}$ and $c^{-}_{j}, d^{-}_{j}$ 
$(j=k+1,\ldots ,m)$ 
are obtained from corresponding generators of $H_{1}(R)$ 
using the half transversal flow associated with $f$,  
see Figure \ref{fig:flow}. 
Let $c'_1,\ldots ,c'_{m}$, 
$d'_1,\ldots, d'_{m}$ 
be the images of 
$c^{-}_1,\ldots ,c^{-}_{m}$, 
$d^{-}_1,\ldots, d^{-}_{m}$ 
in $H_{1}(\partial_{+}W)$.
Then we may describe:
\begin{align*}
&
%\begin{pmatrix}
\Big(
 c'_1 \quad
 \cdots \quad 
 c'_{k} \quad
 c'_{k+1} \quad
 \cdots \quad
 c'_{m} \quad
 d'_1 \quad
 \cdots \quad
 d'_{k} \quad
 d'_{k+1} \quad
 \cdots \quad 
 d'_{m} 
%\end{pmatrix}^T
\Big)^T \\ 
&
=G_1
%\begin{pmatrix}
\Big(
	c^{+}_1 \quad
	\cdots \quad 
	c^{+}_{k} \quad 
	c^{+}_{k+1} \quad
	\cdots \quad
	c^{+}_{m} \quad
	d^{+}_1 \quad
	\cdots \quad
	d^{+}_{k} \quad
	d^{+}_{k+1} \quad
	\cdots \quad 
	d^{+}_{m} 
%\end{pmatrix}^T 
\Big)^T.
\end{align*}

We call $G_1$ the {\it monodromy matrix\/}.
For $n\ge 1$, 
we have: 
$$
\Big(
	g^{n}_{*}(c^{+}_1) 
	\cdots 
	g^{n}_{*}(c^{+}_{k}) \,
	g^{n}_{*}(c^{+}_{k+1}) 
	\cdots 
	g^{n}_{*}(c^{+}_{m}) \,
	g^{n}_{*}(d^{+}_{1}) 
	\cdots 
	g^{n}_{*}(d^{+}_{k}) \,
	g^{n}_{*}(d^{+}_{k+1}) 
	\cdots 
	g^{n}_{*}(d^{+}_{m})\Big)^T
$$
$$
=G_1^{n} 
\Big(
	c^{+}_1 \quad
	\cdots \quad
	c^{+}_{k} \quad
	c^{+}_{k+1} \quad
	\cdots \quad
	c^{+}_{m} \quad
	d^{+}_1 \quad
	\cdots \quad
	d^{+}_{k} \quad
       d^{+}_{k+1} \quad
	\cdots \quad
	d^{+}_{m} 
\Big)^T
$$

Here $(\cdot)^T$ stands for the transposition of a matrix. 

The monodromy $g$ is an orientation preserving diffeomorphism between surfaces, 
then $G_1\in Sp(2m,\mathbb Z)$, in particular 
$\det G_1$=1. 
%Thus the eigenvalues of $F_1$ are the root of unity. 
Further $R$ is a closed or once punctured surface in our setting. 
If $R$ is closed, then ${\rm trace}(G_0)={\rm trace}(G_2)=1$. So, 
if $|t|$ is sufficiently small, 
\begin{align*}
\zeta_{g}(t)
&=\exp\big(\sum_{k=1}^{\infty}\frac{t^k}{k}(2-\text{trace } G_1^k)\big) \\
&=\exp\big(\log(1-t)^{-2}+\text{trace}(\log(I-t\cdot G_1))\big)\\
&=\frac{\det(I-t\cdot G_1)}{(1-t)^2}
\end{align*}
If $R$ is a once punctured surface, we have: 
$$\zeta_{g}(t)=\frac{\det(I-t\cdot G_1)}{1-t}$$ 
by the same argument, if $|t|$ is sufficiently small.
Here $I$ stands for the identity matrix.

%--------------Section7.tex-----'07/09/17------

\section{Counting flow lines}

In this section, we consider counting gradient flow lines from 
critical points of index 2 to those of index 1, 
which are obtained from a circle-valued Morse map $M\rightarrow S^1$, 
according to Section 2.

In our setting, there are only critical points of index 1 and 2, 
we can observe the torsion $\tau_g(t)$ of the chain complex 
((\ref{eqn:chain})
$0\longleftarrow \mathcal N_1\overset{D}{\longleftarrow}\mathcal N_2\longleftarrow 0$)
as follows.

\begin{figure}
\centering
\includegraphics[width=.3\textwidth]{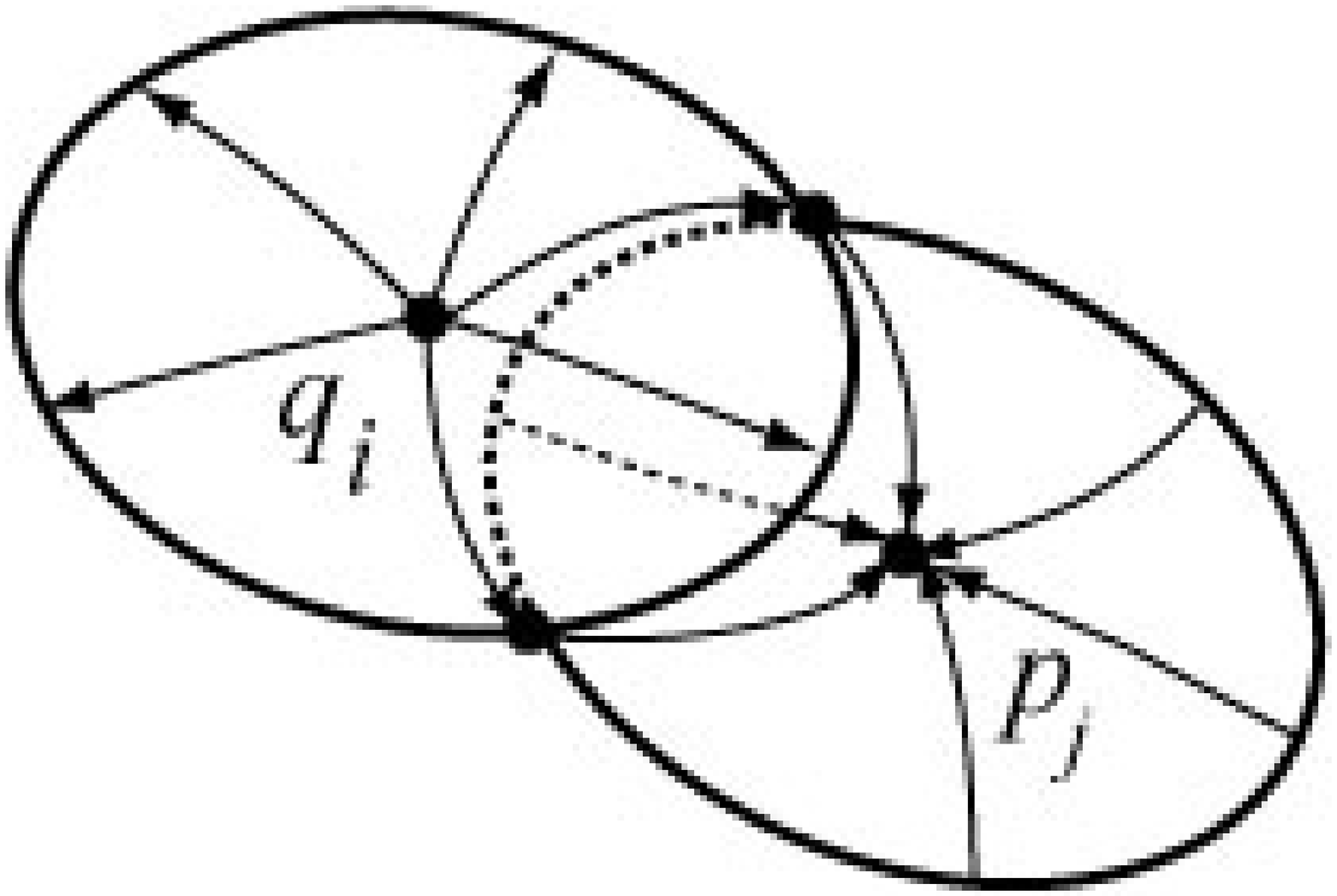}
\caption{}\label{fig:critical}
\end{figure}

As in the previous sections, 
we consider only a monodromy matrix which is obtained from 
a symmetric Heegaard splitting and a half-transversal flow. 
The Novikov module $\mathcal N_{1}$ ($\mathcal N_{2}$ resp.) 
of the pair $(f,v)$ is generated by  $S_{1}(f)=\{p_{1},\ldots,p_{k}\}$ 
($S_{2}(f)=\{q_{1},\ldots,q_{k}\}$ resp.), 
i.e., 
the center points of the disk bounded $c_{i}$ ($tc_{i}$ resp.) 
($i=1,2,\ldots,k$). 
See Figure \ref{fig:compression} and \ref{fig:symmetricflow}.
Therefore the $i\times(m+j)$th-component of the matrix 
$G_{1}$ stands for
the algebraic number of the flow lines between 
$q_{i}$ and $p_{j}\,(1\le i, j\le k)$. 
See Figure \ref{fig:critical} for the schematic image.
Let $D^{(n)}_{ij}$ be the $i\times (m+j)$th-component of $G^{n}_{1},\,
\,(1\le i, j\le k)$.
Then we have:

\begin{definition}\label{def:tau}
We define 
$$\tau_{g}(t)=\det(D_{ij}(t)), 
\hspace{0.3cm}
\text{  where  }
D_{ij}(t)=\sum_{n=1}^{\infty}(D^{(n)}_{ij} \cdot t^{n-1}), \,
1\le i, j\le k.$$
If $M$ has no critical points, 
i.e., $M$ is the fibre bundle over $S^1$ with fibre $R$, then  
$\tau_{g}(t)$ is defined to be 1. 
\end{definition}

%Let $\overline{M}$ be the infinite cyclic covering of $M$ 
%associated a copmact surface $R=h^{-1}(t) (\in H^1(M))$ 
%where $t$ is a regular value of $h$.  
%We denote by $p$ the covering map and by $\overline{R}_i$ 
%the image of $p^{-1}(R)$. 

By taking $|t|$ sufficiently small, 
we have: 
$$\sum^{\infty}_{k=1}G_{1}^{n}\cdot t^{n-1}=G_1(I-t\cdot G_1)^{-1}.$$
Therefore, $D_{ij}(t)$ is the $i\times(m+j)$th-component of 
$G_1(I-t\cdot G_1)^{-1},\,\,(1\le i, j\le k).$

We present the concrete examples for $\tau_g(t)$ in Section \ref{s:examples}.

%--------------Section8.tex-----'07/09/18-----

\section{Examples}
\label{s:examples}

\begin{figure}
\centering
\includegraphics[width=.6\textwidth]{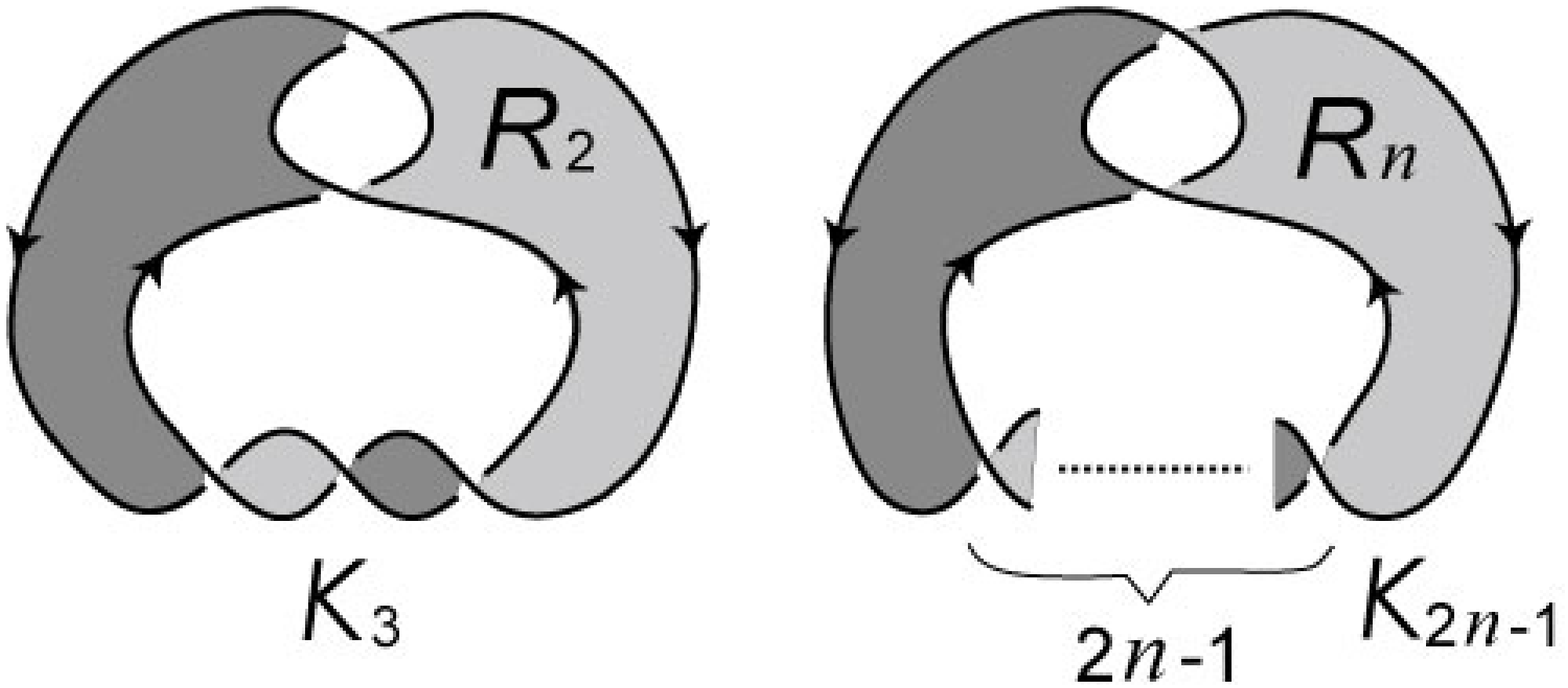}
\caption{}\label{fig:twist}
\end{figure}

In this section, we 
%demonstrate the contents in the previous sections using 
consider twist knots $\mathcal K_{2n-1}\, (n=1,2,3,\ldots)$. 
Note that the Alexander polynomial of $\mathcal K_{2n-1}$ 
is $-n+(2n-1)t-nt^2$.
A twist knot has a genus one Seifert surface $R_n$ as illustrated in 
Figure \ref{fig:twist}.  The twist knot $\mathcal K_1$ is 
the trefoil knot, then it is fibred and treated in Example \ref{ex:trefoil}. 
So, we assume that $n\ge 2$.

Let $X_n$ be the complement of the knot $\mathcal K_{2n-1}$.  

\begin{lemma}\label{lem:example1}
$\mathcal M\mathcal N(X_n,R_n)=2$ for any $n\, (n=2,3,\ldots)$.
\end{lemma}

\begin{proof}
Let $\lambda$ and $\lambda'$ be arcs whose boundaries are in $R_n$ 
as illustrated in Figure \ref{fig:spine}, and $(X_n, R_+, R_-)$ the 
sutured manifold for $R_n$. 
Note that $\partial\lambda=\partial\lambda'$, and
$R_+$ ($R_-$ resp.)  
intersects $\lambda$ ($\lambda'$ resp.) 
transversely in one point. 
Then the regular neighborhood of 
$R_+\cup \lambda$ and $R_-\cup\lambda'$ in $X_n$ are compression bodies. 
Therefore we have only to show that 
the sutured manifold cl$((X_n,R_+,R_-)-(N(R_+\cup\lambda)\cup N(R_-\cup\lambda')))$, 
denoted by  $(\breve{X}_n,\breve{R}_+,\breve{R}_-)$, 
is a product sutured manifold.
We consider the case of $\mathcal K_{5}$ $(n=3)$ since the other cases can be seen 
by the same method.

\begin{figure}
\centering
\includegraphics[width=.6\textwidth]{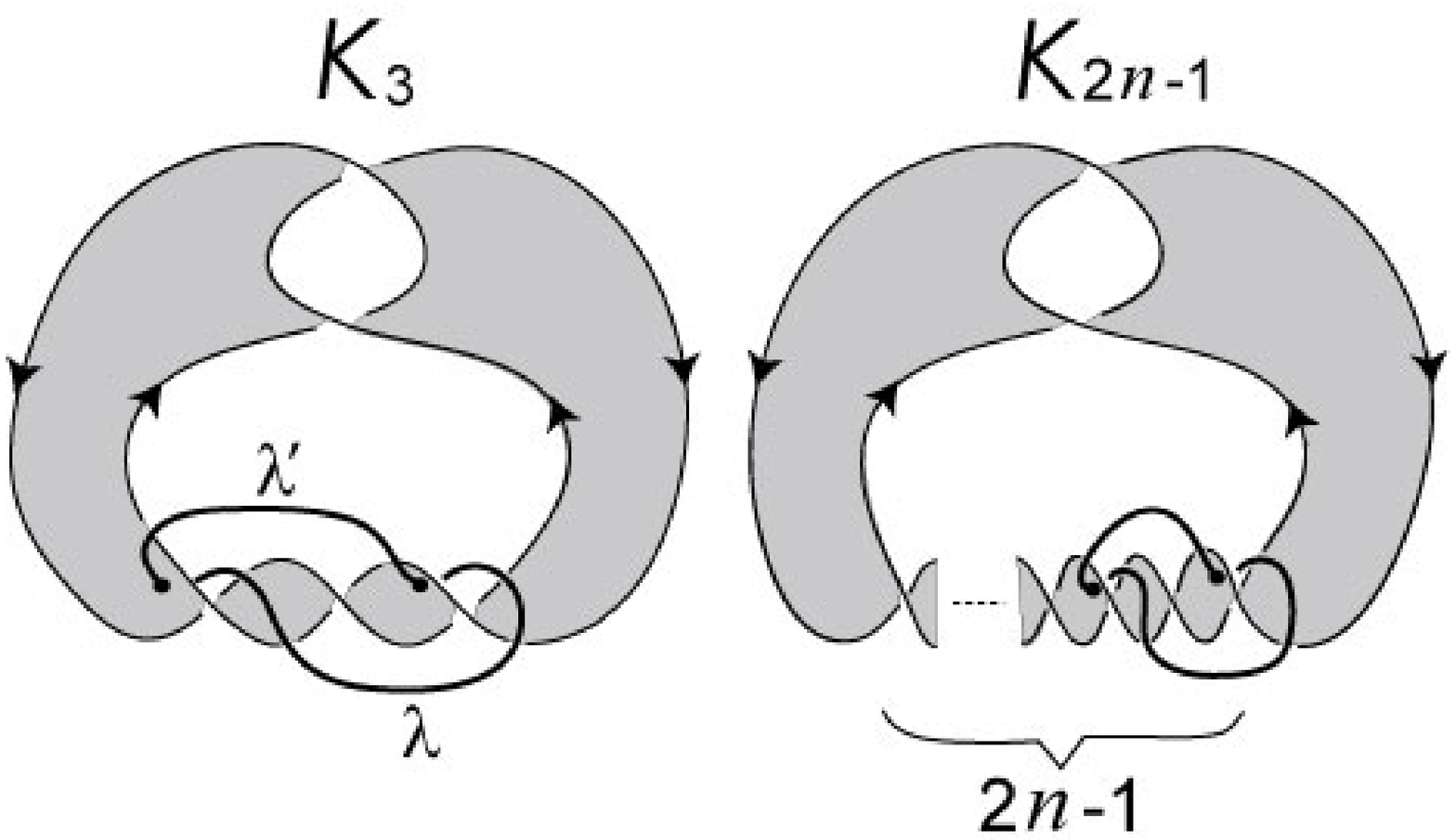}
\caption{}\label{fig:spine}
\end{figure}

\begin{figure}
\centering
\includegraphics[width=.7\textwidth]{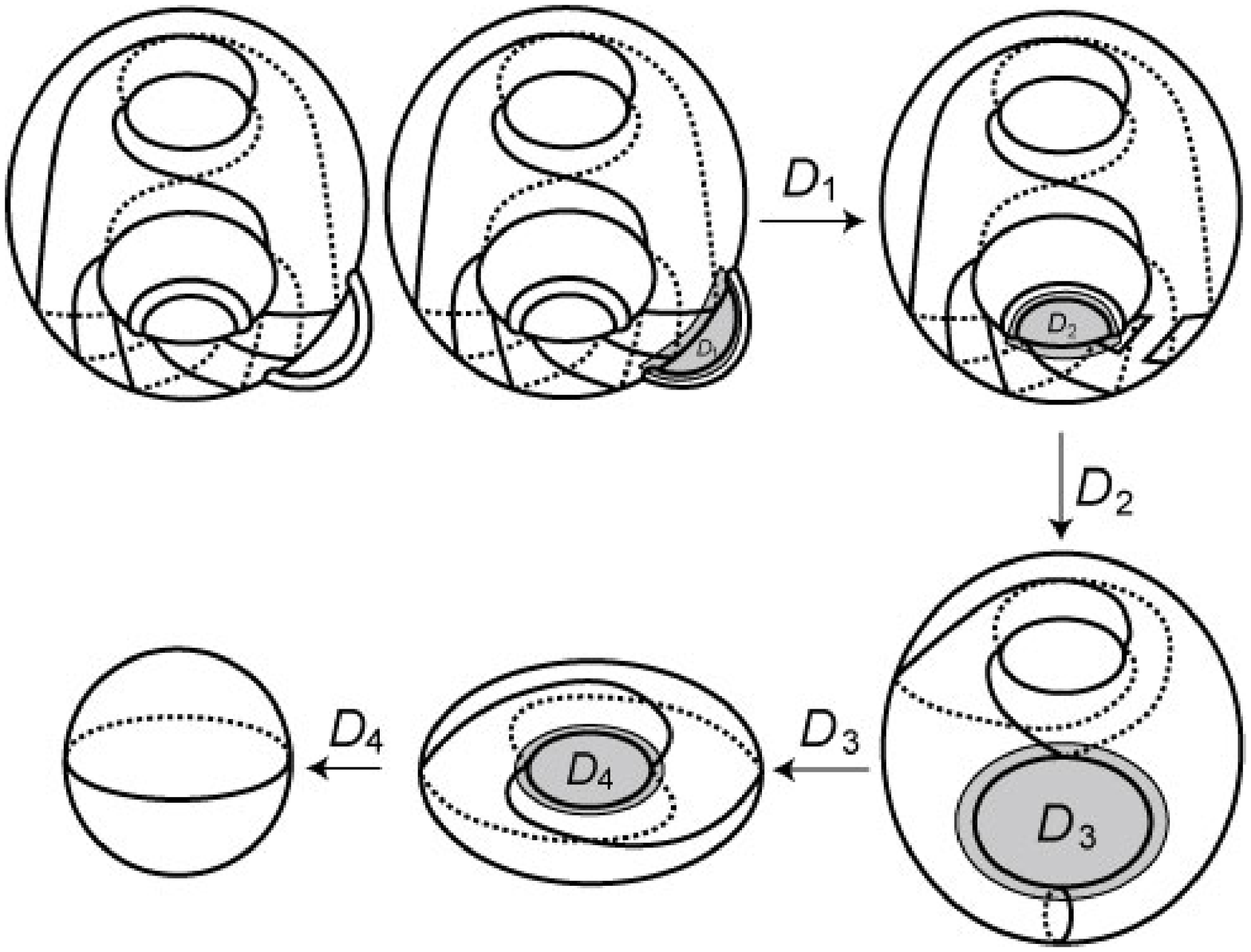}
\caption{}\label{fig:K5}
\end{figure}

\begin{figure}
\centering
\includegraphics[width=.35\textwidth]{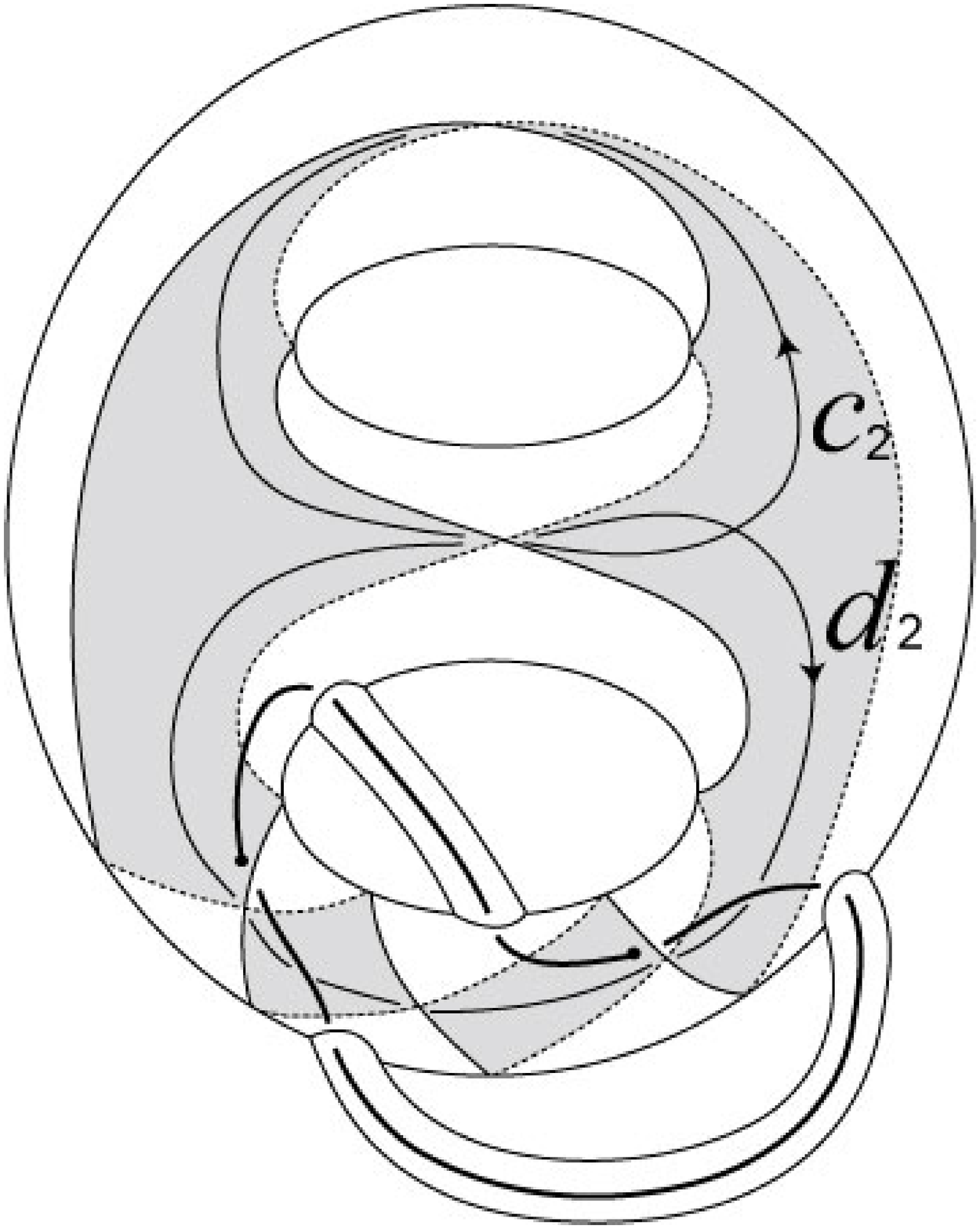}
\caption{}\label{fig:symmetry}
\end{figure}

Let $D_1$ be the product disk in $(\breve{X}_3,\breve{R}_+,\breve{R}_-)$ as 
illustrated in Figure \ref{fig:K5} (shaded part), that is, 
the disk $D_1$ is properly embedded disk in $\breve{X}_3$ 
such that $\partial D_1\cap \breve{R}_+$ ($\partial D_1\cap \breve{R}_- resp.)$ 
is an arc properly embedded in 
$\breve{R}_+$ ($\breve{R}_-$ resp.). 
We decompose $\breve{X}_3$ along $D_1$ and connect the suture naturally, 
then we obtain a new sutured manifold 
$(\breve{X}^1_3,\breve{R}^1_+,\breve{R}^1_-)$.
This decomposition is called a product decomposition \cite{gabai2}. 
Similarly, we decompose $(\breve{X}^1_3,\breve{R}^1_+,\breve{R}^1_-)$ 
along the product disk $D_2$, then we have a sutured manifold 
$(\breve{X}^2_3,\breve{R}^2_+,\breve{R}^2_-)$. 
See Figure \ref{fig:K5}. Thus we have a sequence of the product decompositions: 
$$(\breve{X}_3,\breve{R}_+,\breve{R}_-)
\overset{D_1}{\rightarrow}
(\breve{X}^1_3,\breve{R}^1_+,\breve{R}^1_-)
\overset{D_2}{\rightarrow}
(\breve{X}^2_3,\breve{R}^2_+,\breve{R}^2_-)
\overset{D_3}{\rightarrow}
(\breve{X}^3_3,\breve{R}^3_+,\breve{R}^3_-)
\overset{D_4}{\rightarrow}
(\breve{X}^4_3,\breve{R}^4_+,\breve{R}^4_-),$$
where 
$\breve{X}^4_3$ is homeomorphic to the 3-ball and 
both $\breve{R}^4_+$ and $\breve{R}^4_-$ are disks. 
This shows that 
$(\breve{X}_3,\breve{R}_+,\breve{R}_-)$ is a product sutured manifold 
by \cite{gabai2}. 
By the same argument, we have that 
$(\breve{X}_n,\breve{R}_+,\breve{R}_-)$ is a product sutured manifold.
This completes the proof.
\end{proof}

We denote by $(W_n,W'_n)$
the Heegaard splitting of $(X_n,R_{+},R_{-})$, 
which is obtained in the proof of Lemma \ref{lem:example1}.

\begin{lemma}\label{lem:example2}
The Heegaard splitting $(W_n,W'_n)$ is symmetric.
\end{lemma}

\begin{proof}
Since $\partial\lambda=\partial\lambda'$ 
and $(\breve{X}_{n},\breve{R}_{+},\breve{R}_{-})$ 
is a product sutured manifold, we have this lemma.
\end{proof}

For the simplicity, we discuss the case of $\mathcal K_{3}$ 
$(n=2)$ 
in the next lemma. 
The general case can be obtained by the same method. 

\begin{lemma}
The Heegaard splitting $(W_2,W'_2)$ induces a monodromy matrix 
presented by
$$
G_1=
\begin{pmatrix}
\hphantom{-}1 & \hphantom{-}1 & -2 & -1 \\
\hphantom{-}0 & \hphantom{-}1 & -1 & \hphantom{-}0 \\
\hphantom{-}0 & \hphantom{-}0 & \hphantom{-}1 & \hphantom{-}0\\
\hphantom{-}0 & -1 & \hphantom{-}0 & \hphantom{-}1 
\end{pmatrix}.
$$
Moreover, 
we have 
$$\zeta_{g}(t)=(1-t)^3 \text{ and } 
\tau_{g}(t)=\frac{-2+3t-2t^{2}}{(1-t)^4}.$$
\end{lemma}

\begin{figure}
\centering
\includegraphics[width=.7\textwidth]{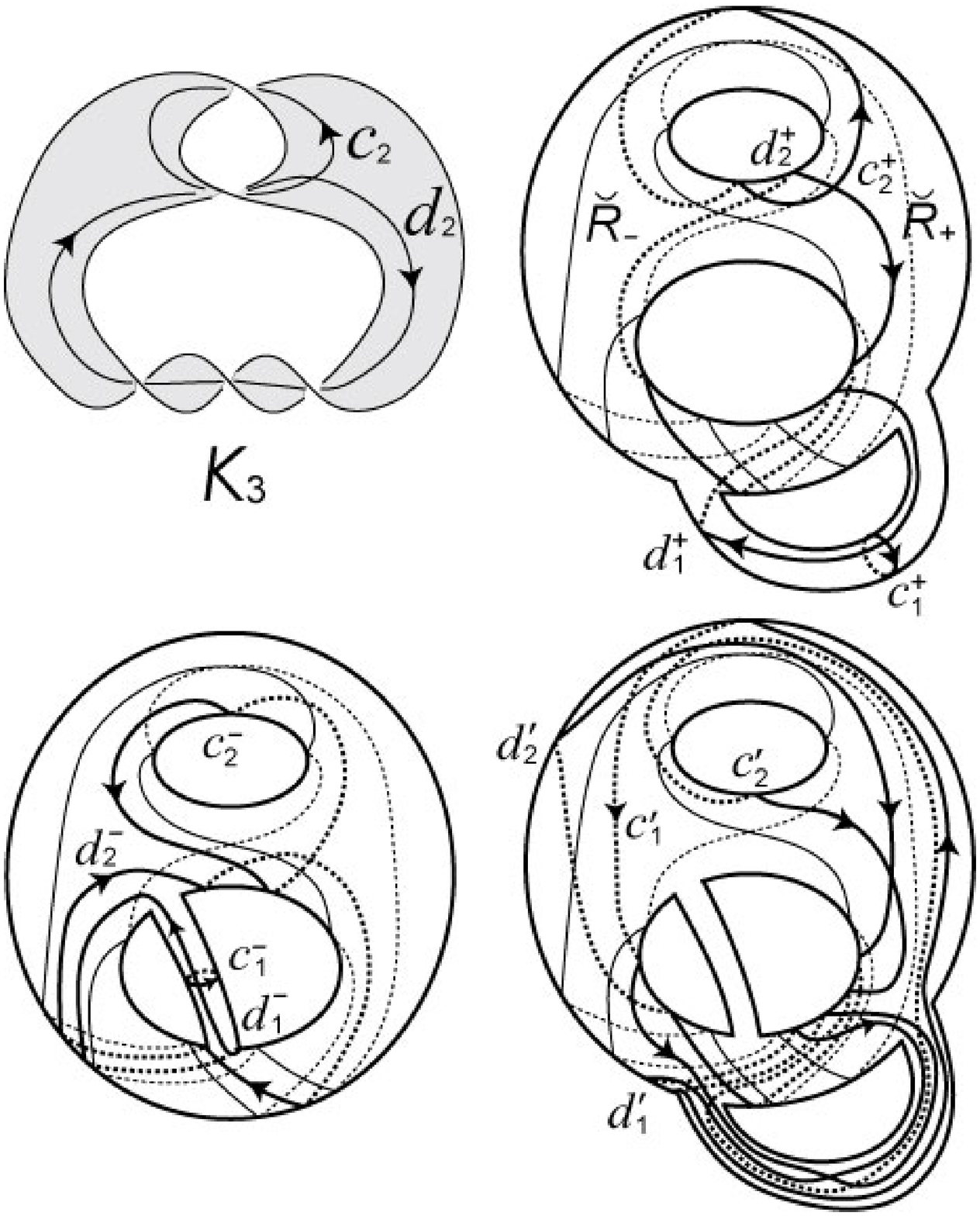}
\caption{}\label{fig:generator}
\end{figure}

\begin{proof}
We take a basis $c_{2},\,d_{2}$ of $H_1(R)$ 
as illustrated in Figure \ref{fig:symmetry}, 
then we have a basis $c^{+}_{2},\,d^{+}_{2}$ of $H_1(R_+)$ 
($c^{-}_{2},\,d^{-}_{2}$ of $H_1(R_-)$ resp.) 
as in the upper right-hand figure 
(lower left-hand figure resp.) in Figure \ref{fig:generator}. 
Note the positions of $\lambda, \lambda'$ and $c_{2},\,d_{2}$ in Figure \ref{fig:symmetry}. 
Let $(\breve{X}_2,\breve{R}_+,\breve{R}_-)$ be the sutured manifold 
cl$(X_2,R_+,R_-)-(N(R_+\cup\lambda)\cup N(R_-\cup\lambda'))$ 
as in the proof of Lemma \ref{lem:example1}.
Here we see that $c^{+}_{1},\,c^{+}_{2},\,d^{+}_{1},\,d^{+}_{2}\subset \breve{R}_+$ 
and $c^{-}_{1},\,c^{-}_{2},\,d^{-}_{1},\,d^{-}_{2}\subset \breve{R}_-$.  
Since $(\breve{X}_2,\breve{R}_+,\breve{R}_-)$ is a product sutured manifold, 
we can move $c^{-}_{1},\,c^{-}_{2},\,d^{-}_{1},\,d^{-}_{2}$ by a free homotopy 
from $\breve{R}_{-}$ to $\breve{R}_{+}$. 
We denote their images by 
$c'_{1},\,c'_{2},\,d'_{1},\,d'_{2}$. 
Then we can see that 
they sit as in the lower right-hand figure in Figure \ref{fig:generator}.
Hence we have:
$ c'_{1}=c^{+}_{1}+c^{+}_{2}-2d^{+}_{1}-d^{+}_{2},\,
  c'_{2}=c^{+}_{2}-d^{+}_{1},\,
  d'_{1}=d^{+}_{1},\,
  d'_{2}=-c^{+}_{2}+d^{+}_{2}.$
Therefore we have the monodromy matrix $G_{1}$ 
in the statement of this lemma, 
and we have
$$
\zeta_{g}(t)
=\frac{\det(I-t\cdot G_1)}{1-t}
=\frac{(1-t)^4}{1-t}
=(1-t)^3.
$$
Note that the convergence radius is 1. 

On the other hand, 
$$
G_1(I-t\cdot G_1)^{-1}
=
\begin{pmatrix}
\frac{1}{(1-t)} & \frac{1}{(1-t)^{3}} & \frac{-2+3t-2t^2}{(1-t)^4} & \frac{-1}{(1-t)^2} \\
0 & \frac{1}{1-t} & \frac{-1}{(1-t)^2} & 0 \\
0 & 0 & \frac{1}{1-t} & 0 \\
0 & \frac{-1}{(1-t)^2} & \frac{t}{(1-t)^3} & \frac{1}{1-t} 
\end{pmatrix}.
$$
Thus we have 
$\displaystyle{
\tau_{g}(t)=
\frac{-2+3t-2t^2}{(1-t)^4}}.$
\end{proof}

By the same argument, we have: 

\begin{proposition}
The Heegaard splitting $(W_n,W'_n)$ induces a monodromy matrix 
presented by
$$
\begin{pmatrix}
\hphantom{-}1 & \hphantom{-}1 & -n & -1 \\
\hphantom{-}0 & \hphantom{-}1 & -1 & \hphantom{-}0  \\
\hphantom{-}0 & \hphantom{-}0 & \hphantom{-}1 & \hphantom{-}0 \\
\hphantom{-}0  & -1 & \hphantom{-}0 & \hphantom{-}1 
\end{pmatrix}.
$$
Moreover, 
we have 
$$\zeta_{g}(t)=(1-t)^3 \text{ and } 
\tau_{g}(t)=\frac{-n+(2n-1)t-nt^{2}}{(1-t)^4}.$$
\end{proposition}

\begin{figure}
\centering
\includegraphics[width=.35\textwidth]{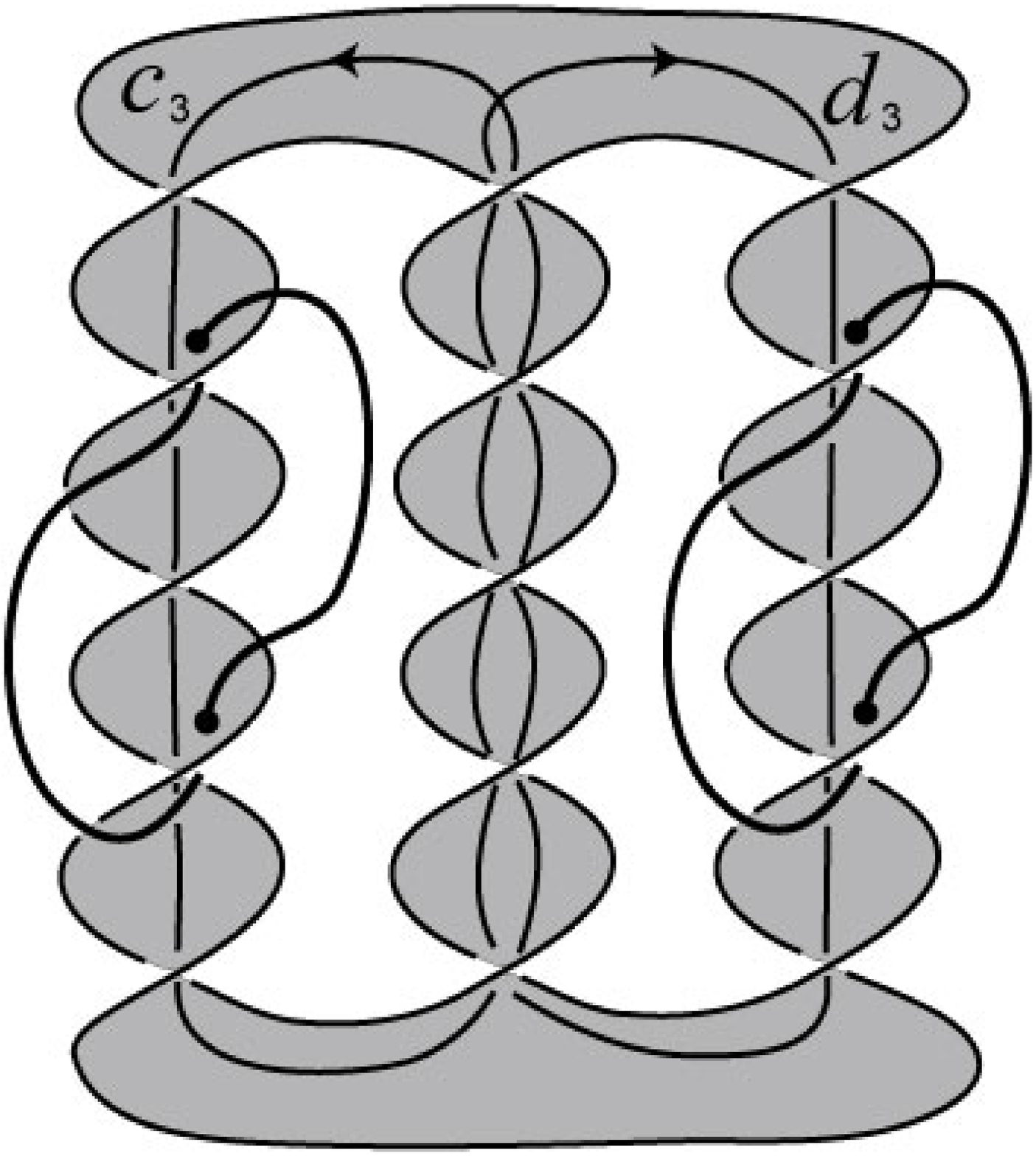}
\caption{}\label{fig:pretzel}
\end{figure}

\begin{example}
Let $K$ be the pretzel knot of type $(5,5,5)$ and 
we consider the symmetric Heegaard splitting associated with 
Figure \ref{fig:pretzel}. Then, 
$$G_1=
\begin{pmatrix}
1 & 0 & 1 & -5 & -2 & 0 \\
0 & 1 & 0 & -3 & -5 & 1 \\
0 & 0 & 1 & \hphantom{-}0 & \hphantom{-}1 & 0 \\
0 & 0 & 0 & \hphantom{-}1 & \hphantom{-}0 & 0 \\
0 & 0 & 0 & \hphantom{-}0 & \hphantom{-}1 & 0 \\
0 & 0 & 0 & -1 & \hphantom{-}0 & 1 \
\end{pmatrix}.
$$
Thus we have $\zeta_g(t)=(1-t)^5$. 
Further, 
$$D_{11}(t)=D_{22}(t)=\frac{-5}{(1-t)^2},\, 
  D_{12}(t)=\frac{-2+3t}{(1-t)^3},\,
  D_{21}(t)=\frac{-3+2t}{(1-t)^3},\,$$
$$\tau_g(t)=\frac{19-37t+19t^2}{(1-t)^6}.$$ 
\end{example}

\medskip 

{\bf Acknowledgment.}
A part of this work was carried out while the first and third authors 
were visiting Institut des Hautes \'Etudes Scientifiques in France 
and the third author was visiting Tokyo University 
of Agriculture and Technology supported by Japan Society for the Promotion Science (JSPS). 
Further, the final version was done while the first author visited 
Universit\'e de Nantes and the second author visited Columbia University 
supported by JSPS.
They would like to thank them for their hospitality and support.

\medskip

%--------------ref.tex-----'07/09/11-----

\bibliographystyle{amsplain}

\end{document}

%% file: auxxx.tex
\renewcommand{\a}{\alpha}
\renewcommand{\b}{\beta}
\newcommand{\g}{\gamma}

\newcommand{\z}{\zeta}

\renewcommand{\l}{\lambda}

\newcommand{\s}{\sigma}

\newcommand{\G}{\Gamma}
\newcommand{\D}{\Delta}

\renewcommand{\L}{\Lambda}

\newcommand{\NN}{{\mathcal N}}

\newcommand{\RR}{{\mathcal R}}

\newcommand{\WW}{{\mathcal W}}

\newcommand{\ZZ}{{\mathcal Z}}

\newcommand{\FFF}{{\mathbf{F}} }

\newcommand{\NNN}{{\mathbf{N}}}

\newcommand{\QQQ}{{\mathbf{Q}}}
\newcommand{\RRR}{{\mathbf{R}}}

\newcommand{\ZZZ}{{\mathbf{Z}}}

\newcommand{\Wh}{\text{\rm Wh}}

\newcommand{\ind}{\text{\rm ind\hspace{0.05cm}}}

\newcommand{\Int}{\text{\rm Int }}

\newcommand{\bere}{\begin{rema}}
\newcommand{\bede}{\begin{defi}}

\renewcommand{\beth}{\begin{theo}}
\newcommand{\bele}{\begin{lemm}}
\newcommand{\bepr}{\begin{prop}}
\newcommand{\beeq}{\begin{equation}}
\newcommand{\bega}{\begin{gather}}
\newcommand{\begaa}{\begin{gather*}}
\newcommand{\been}{\begin{enumerate}}

\newcommand{\bedee}{\begin{defii}}
\newcommand{\bethh}{\begin{theoo}}
\newcommand{\belee}{\begin{lemmm}}
\newcommand{\beprr}{\begin{propp}}

\newcommand{\beco}{\begin{coro}}

\newcommand{\beal}{\begin{aligned}}

\newcommand{\enre}{\end{rema}}

\newcommand{\enco}{\end{coro}}
\newcommand{\enpr}{\end{prop}}
\newcommand{\enth}{\end{theo}}
\newcommand{\enle}{\end{lemm}}
\newcommand{\enen}{\end{enumerate}}
\newcommand{\enga}{\end{gather}}
\newcommand{\engaa}{\end{gather*}}
\newcommand{\eneq}{\end{equation}}
\newcommand{\enal}{\end{aligned}}

\newcommand{\bq}{\begin{equation}}
\newcommand{\bqq}{\begin{equation*}}

\newcommand{\mx}{\mbox}

\newcommand{\wi}{\widetilde}

\newcommand{\ove}{\overline}

\newcommand{\wh}{\widehat}

\newcommand{\sm}{\setminus}

\newcommand{\sbs}{\subset}

\newcommand{\tens}[1]{\underset{#1}{\otimes}}

\newcommand{\Prf}{{\it Proof.\quad}}

\newcommand{\chart}{\Phi_p:U_p\to B^n(0,r_p)}
\newcommand{\atlas}{\{\Phi_p:U_p\to B^n(0,r_p)\}_{p\in S(f)}}

\newcommand{\pr}{\partial}

\newcommand{\qs}{\hfill\square}